\documentclass[aop,MSNbibl,seceqn,citesort,dvips]{arximspdf}
\usepackage{graphicx}

\doi{10.1214/11-AOP713} 
\volume{41}
\issue{2}
\pubyear{2013}
\firstpage{636}
\lastpage{670}

\makeatletter

\newtheorem{theorem}{Theorem}[section]
\newtheorem{lemma}[theorem]{Lemma}
\newtheorem{proposition}[theorem]{Proposition}
\newtheorem{corollary}[theorem]{Corollary}

\newproclaim{assumption}[theorem]{Assumption}

\newcommand{\E}{\mathbf{E}}
\newcommand{\N}{\mathbf{N}}
\newcommand{\Z}{\mathbf{Z}}
\newcommand{\p}{\mathbf{P}}
\newcommand{\R}{\mathbf{R}}
\newcommand{\T}{\mathbf{T}}

\newcommand{\CA}{\mathcal{A}}
\newcommand{\CB}{\mathcal{B}}
\newcommand{\CR}{\mathcal{R}}
\newcommand{\CS}{\mathcal{S}}
\newcommand{\CT}{\mathcal{T}}

\newcommand{\var}{\operatorname{Var}}
\newcommand{\cov}{{\mathrm{cov}}}

\newcommand{\one}{\mathbf{1}}

\newcommand{\ol}{\overline}
\newcommand{\wt}{\widetilde}

\newcommand{\mix}{{\mathrm{mix}}}

\makeatother

\begin{document}
\begin{frontmatter}

\title{Painting a graph with competing random walks\thanksref{T2}}
\runtitle{Painting a graph with competing random walks}

\begin{aug}
\author[A]{\fnms{Jason} \snm{Miller}\corref{}\ead[label=e2]{jmiller@math.stanford.edu}}
\runauthor{J. Miller}
\affiliation{Stanford University}
\address[A]{Department of Mathematics\\
Stanford University\\
Stanford, California 94305\\
USA\\
\printead{e2}} 
\end{aug}

\thankstext{T2}{Supported in part by NSF
Grants DMS-04-06042 and DMS-08-06211.}

\received{\smonth{3} \syear{2010}}
\revised{\smonth{4} \syear{2011}}

%
\begin{abstract}
Let $X_1, X_2$ be\vspace*{1pt} independent random walks on $\Z_n^d$, $d
\geq3$, each starting from the uniform distribution. Initially, each
site of $\Z_n^d$ is unmarked, and, whenever $X_i$ visits such a site,
it is set irreversibly to $i$. The mean of $|\CA_i|$, the cardinality
of the set $\CA_i$ of sites painted by $i$, once all of $\Z_n^d$ has
been visited, is $\frac{1}{2}n^d$ by symmetry. We prove the following
conjecture due to Pemantle and Peres: for each $d \geq3$ there exists a
constant $\alpha_d$ such that $\lim_{n \to\infty} \var(|\CA_i|) /
h_d(n) = \frac{1}{4}\alpha_d$ where $h_3(n) = n^4$, $h_4(n) = n^4 (\log
n)$ and $h_d(n) = n^d$ for $d \geq5$. We will also identify $\alpha_d$
explicitly and show that $\alpha_d \to1$ as $d \to\infty$. This is a
special case of a more general theorem which gives the asymptotics of
$\var(|\CA_i|)$ for a large class of transient, vertex transitive
graphs; other examples include the hypercube and the Caley graph of the
symmetric group generated by transpositions.
\end{abstract}

%

\end{frontmatter}

\section{Introduction}\label{intro}

Suppose that $X_1,X_2$ are independent random walks on a graph $G =
(V,E)$ starting from stationarity. Initially, each vertex of $G$ is
unmarked, and, whenever $X_i$ visits such a site, it is marked $i$
irreversibly. If both $X_1$ and $X_2$ visit a site for the first time
simultaneously, then the mark is chosen by the flip of an independent
fair coin. Let $\CA_i$ be the set of sites marked $i$ once every vertex
of $G$ has been visited. By symmetry, it is obvious that $\E|\CA_i| =
\frac{1}{2} |V|$. The purpose of this manuscript is to derive precise
asymptotics for $\var(|\CA_i|)$ for many families of graphs.

The process by which a single random walk covers a graph has been
studied extensively. Examples of interesting statistics include the
expected amount of time it takes for the random walk to visit every
site~\cite{M88,DPRZCOV04}, the growth exponent of the set of sites
visited most frequently~\cite{DPRZTHICK01} and the clustering and
correlation structure of the last visited points
\cite{BH91,DPRZLATE06,MP11}. The motivation for this work is to understand
better how multiple random walks cover a graph.

The investigation of the statistical properties of $\CA_i$ was first
proposed in the work of Gomes Jr. et al.~\cite{JLSH96}. Their
motivation was to study the technical challenges associated with
physical problems involving interacting random walks. They estimate the
growth of $\E|\CB|$ where $\CB$ is the interface separating $\CA_1$
from $\CA_2$ in the special case of the one-cycle $\Z_n^1$. As with
$\E
|\CA_i|$, computing $\E|\CB|$ for $\Z_n^d$ becomes trivial for $d
\geq
3$ since it is easy to see that, with probability strictly between $0$
and~$1$, for any pair of adjacent vertices $x,y$, $X_2$ will hit $y$
before $X_1$, conditional on the event that $X_1$ hits $x$ first. On
the other hand, estimating $\var(|\CA_i|)$ in this setting is
challenging since its expansion in terms of correlation functions
exhibits significant cancellation which, when ignored, leads to bounds
that are quite imprecise. We will develop this point further at the end
of the \hyperref[intro]{Introduction}.

The problem we consider here was formulated by Hilhorst, though in a
slightly different setting. Rather\vspace*{1pt} than considering the
sets of sites $\CA_1, \CA_2$ \textit{first painted} by $X_1,X_2$,
respectively, it is also natural to study the sets $\wt{\CA}_1,
\wt{\CA}_2$ of sites \textit{most recently painted} by $X_1, X_2$,
respectively. In other words, in the latter formulation the constraint
that the marks are irreversible is removed. It turns\vspace*{1pt} out
that these two classes of problems are equivalent, which is to say
$(\CA_1, \CA_2) \stackrel{d}{=} (\wt {\CA }_1, \wt{\CA}_2)$. This
helpful observation, which follows from the time-reversibility of
random walk, was made and communicated to us by Comets.

We restrict our attention to \textit{lazy} walks $X_1,X_2$ to avoid
issues of periodicity, and in particular to ensure that the random walk
has a unique stationary distribution. That is, the one-step transition
kernel is given by
\[
p(x,y;G) =
\cases{
\dfrac{1}{2}, &\quad if $x = y$,\cr
\dfrac{1}{2\deg(x)}, &\quad if $x \sim y$,\vspace*{2pt}\cr
0, &\quad otherwise,}
\]
where $x \sim y$ means that $x$ is adjacent to $y$ in $G$. The
particular choice of holding probability $\frac{1}{2}$ is not
important for the proof; indeed, any $\lambda\in(0,1)$ would suffice.
Our proofs also work in the setting of continuous time walks. Let
$p^t(\cdot,\cdot;G)$ be the $t$-step transition kernel of a lazy random
walk on $G$ and $\pi(\cdot;G)$ its unique stationary distribution.

Our main result is the precise asymptotics for $\var(|\CA_i|)$ on tori
of dimension at least three, thus verifying a conjecture due to
Pemantle and Peres~\cite{DICKERTHESIS}, page 35.
%
%
\begin{theorem}
\label{thmtori}
Suppose that $G_n = \Z_n^d$, $d \geq3$. There exists a finite constant
$\alpha_d > 0$ such that
\[
\lim_{n \to\infty} \frac{\var(|\CA_i|)}{h_d(n)} = \frac
{1}{4}\alpha_d,
\]
where
\[
h_d(n) = \cases{
n^4, &\quad if $d =3$,\cr
n^4 (\log n), &\quad if $d = 4$, and\cr
n^d, &\quad if $d \geq5$.}
\]
\end{theorem}

Our proof allows us to identify $\alpha_d$ explicitly and is given as
follows. Let
%
%
\begin{equation}
\label{eqngreensfunctionzd}
G(x;\Z^d) = \E_0 \sum_{t=0}^\infty\one_{\{X(t) = x\}}
\end{equation}
be the Green's function for lazy random walk on $\Z^d$. This is the
amount of time a random walk initialized at $0$ spends at $x$ before
escaping to $\infty$. For $d \geq5$,
%
%
\begin{equation}
\label{eqnalphad5}
\alpha_{d} = \frac{1}{G^2(0;\Z^d)}\sum_{y \in\Z^d} G^2(y;\Z^d).
\end{equation}
It is not difficult to see that $\alpha_d \to1$ as $d \to\infty$, so
that $\var(|\CA_1|) \approx\frac{1}{4} n^d$ for $d$ and $n$ large is
close to the variance of an i.i.d. marking. For $d = 4$,
%
%
\begin{equation}
\label{eqnalphad4}
\alpha_4 = \lim_{n \to\infty} \frac{1}{G^2(0;\Z^4) \log n} \sum
_{y \in
\Z^4\dvtx|y| \leq n} G^2(y;\Z^4);
\end{equation}
we will explain why this limit exists and is positive and finite in
Proposition~\ref{proplimitexistence}. The definition of $\alpha_3$
is slightly more involved. Let $\T^3$ denote the three-dimensional
continuum torus, $p^t(\cdot,\cdot;\T^3)$ the transition kernel for Brownian motion
on $\T^3$ and
\[
g^T(x,y;\T^3) = \int_0^{T} p^t(x,y;\T^3) \,dt.
\]
Now set
%
%
\begin{eqnarray}
\label{eqnalphad3}
\alpha_3^T &=& \frac{1}{G^2(0;\Z^3)}\int_{\T^3} \int_{\T^3}
\biggl(g^{T/2}(x,y;\T^3) - \frac{1}{2}T\biggr)^2 \,dx
\,dy,\nonumber\\[-8pt]\\[-8pt]
\alpha_3 &=& \lim_{T \to\infty} \alpha_3^T.\nonumber
\end{eqnarray}
The reason that the limit exists and is positive and finite is that
$p^t(x,y;\T^3)$ converges to the uniform density exponentially fast in
$t$; see Proposition~\ref{propmixingdecay} for a discrete version of
this statement.

Throughout the rest of the article, for functions $f,g$, we say that $f
= O(g)$ if there exists constants $c_1, c_2$ such that $|f| \leq c_1 +
c_2 |g|$. We say that $f = \Omega(g)$ if there exists constants $c_1,
c_2$ so that $|f| \geq c_1 + c_2 |g|$. We say that $f = \Theta(g)$ if
$f = O(g)$ and $f = \Omega(g)$. Finally, we say $f = o(g)$ if $\lim_{n
\to\infty} f(n) / g(n) = 0$.

We note that the problem for $d=1$ is trivial: $\var(|\CA_1|) =
\Theta
(n^2)$. Indeed, observe that with positive probability, the distance
between $X_1$ and $X_2$ at time $0$ is at least $\frac{1}{4} n$. In
$cn^2$ steps (for $c$ large enough), $X_1$ has positive probability of
covering the entire cycle while $X_2$ has positive probability of not
leaving an interval of length $\frac{1}{4} n$ containing its starting
point. On this event, $|\CA_1| \geq\frac{3}{4} n$. This proves our
claim as the upper bound is trivial. For $d=2$, the asymptotics of
$\var
(|\CA_1|)$ remains open.

One interesting remark is that the variance for $d=3,4$ is
significantly higher than that of an i.i.d. marking. The results of
Theorem~\ref{thmtori} should also be contrasted with the behavior of
the variance of the range $\CR$ of random walk on $\Z^d$ run up to the
\textit{cover time} $T_\cov(\Z_n^d)$ of $\Z_n^d$, which is the expected
amount of time it takes for a single random walk to visit every site.
When $d \geq3$, $T_\cov(\Z_n^d) \sim c_d n^d (\log n)$; see
\cite{LPW08}. For $d \geq5$, it follows from work of Jain and Orey
\cite{JO68} that $\var(|\CR|) = \Theta( n^d (\log n))$. For $d =3,4$, it
follows from work of Jain and Pruitt~\cite{JP70} that $\var(|\CR|)$ is
$\Theta(n^3 (\log n)^2)$ and $\Theta(n^4 (\log n))$, respectively.

This work opens the doors to many other problems involving two random
walks. Natural next steps include CLTs for the fluctuations of $|\CA
_i|$ and for the number of sites painted by $i$ at time $t$, as well as
the development of an understanding of the geometrical properties of
the clusters of $\CA_i$. The latter seem to be connected to the theory
of \textit{random interlacements}. This is a model developed by Sznitman
in~\cite{S10} to describe the microscopic structure of the points
visited by a random walk on $\Z_n^d$, $d \geq3$, at times $u n^d$ for
$u > 0$---that is, when a constant order of vertices have been
visited. Roughly speaking, the model is a Poisson process on $W^*
\times(0,\infty)$, where $W^*$ is the space of doubly-infinite paths
on $\Z^d$ modulo time-shifts. For a point $(X,U)$ realized in this
process, one should think of $X$ as describing a random walk trajectory
(an ``interlacement'') and $U$ a time parameter. The model was first
developed to study the process of disconnection of a discrete cylinder
by random walk~\cite{DS06} and has been subsequently applied to
understand the fine geometrical structure of random walk in many
different settings~\cite{W08,W10}. Sznitman's theory
generalizes to the setting of $k$ random walks by labeling each
interlacement with an element of $\{1,\ldots,k\}$ i.i.d. at random.
Studying the structure of the clusters in the $\CA_i$ using this
general theory is an interesting research direction.

Theorem~\ref{thmtori} is a special case of a much more general
result, which gives the asymptotics of $\var(|\CA_i|)$ for many other
graphs, such as the hypercube and the Caley graph of the symmetric
group generated by transpositions. We will now review some additional
terminology which is necessary to give a precise statement of the
result. Recall that the \textit{uniform mixing time} of random walk on
$G$ is
\[
T_\mix(G) = \min\biggl\{ t \geq0\dvtx\max_{x,y} \biggl| \frac
{p^t(x,y;G)}{\pi(y;G)} - 1 \biggr| \leq\frac{1}{4} \biggr\},
\]
and the \textit{Green's function} for $G$ is
\[
g(x,y;G) = \sum_{t=0}^{T_\mix(G)} p^t(x,y;G),
\]
that is, the expected amount of time $X_i$ spends at $y$ up until
$T_\mix(G)$ when started from $x$. Let $\tau_i(x) = \min\{t \geq0\dvtx
X_i(t) = x\}$ be the first time $X_i$ hits $x$; we will omit $i$ if
there is only one random walk under consideration. Throughout the rest
of the article, $a \wedge b = \min(a,b)$ for $a,b \in\R$.
%
%
\begin{assumption}
\label{assumpmain}
$(G_n)$ is a sequence of vertex transitive graphs with $|V_n| \to
\infty
$ such that:
{\renewcommand\thelonglist{(\arabic{longlist})}
\renewcommand\labellonglist{\thelonglist}
\begin{longlist}
\item \label{assumpmainmixingbound} $T_\mix(G_n) = o(|V_n| /
({\log}|V_n|)^2)$ and $\lim_{n \to\infty} T_\mix(G_n) = \infty$;
\item \label{assumpmaingreensquarebound} $\sum_{y \neq x_0}
g^2(x_0,y;G_n) = o(T_\mix(G_n) / {\log}|V_n|)$ for each $x_0 \in V_n$ fixed;
\item \label{assumpmainhitmiss} there exists $\rho_0 < 1$ so
that $\p_x[\tau(y) \wedge\tau(z) \leq T_\mix(G_n)] \leq\rho_0$
uniformly in $n$ and $x,y,z \in V_n$ distinct.
\end{longlist}}
\end{assumption}

The purpose of~\ref{assumpmainmixingbound} is that in many cases
we will perform union bounds over time-scales whose length is
proportional to $T_\mix(G_n)$, and the hypothesis gives us explicit
control on how the number of terms in these bounds relates to the size
of~$V_n$. Part~\ref{assumpmaingreensquarebound} gives us
control on the tail behavior of $g$ and, finally, part \ref
{assumpmainhitmiss} says that with uniformly positive probability
the walks we consider do not hit adjacent points within the mixing
time. Note that vertex transitivity implies $g$ is constant along the
diagonal. Part~\ref{assumpmainhitmiss} implies that the number
of times random walk started at $x$ returns to $x$ before the mixing
time is stochastically dominated by a geometric random variable whose
parameter depends only on $\rho_0$. Consequently, we see that there
exists $g_0 > 0$ such that $g(x,x;G_n) \leq g_0$ uniformly in $x$ and $n$.

Assume that $(G_n)$ is a sequence of vertex transitive graphs, and let
%
%
\begin{eqnarray}
\label{eqnfncdef}
f_{n,c}(x,y) &=& \p_x[\tau(y) \leq c T_\mix(G_n)],
\\
\label{eqnfncbardef}
\ol{f}_{n,c} &=& \sum_{y} f_{n,c}(x,y) \pi(y;G_n).
\end{eqnarray}
Note that $\ol{f}_{n,c}$ does not depend on the choice of $x$ since if
we replaced $x$ with $x'$, by vertex transitivity we may precompose
$f_{n,c}$ with an automorphism of $G_n$ which sends $x$ to $x'$.

The general theorem is:
%
%
\begin{theorem}
\label{thmmain}
Suppose that $(G_n)$ satisfies Assumption~\ref{assumpmain}.
Let
\[
F_{n,c} = \sum_{x,y} \bigl(f_{n,c}(x,y) - \ol{f}_{n,c}\bigr)^2.
\]
There exists $\gamma> 0$ so that for every $c \geq2$, we have
%
%
\begin{equation}\label{eqnvariance}
\var( |\CA_i| )
= \bigl(\tfrac{1}{4}+O(\Delta_n)\bigr) F_{n,c}
+ O(e^{-\gamma c} (T_\mix(G_n))^2)
\end{equation}
as $n \to\infty$ where
\[
\Delta_n = \frac{{T_\mix(G_n) \log}|V_n|}{|V_n|}.
\]
\end{theorem}

Applying this to the special cases of the hypercube and the Caley graph
of $S_n$ generated by transpositions leads to the following corollary.
%
%
\begin{corollary}
\label{corhypsym}
Suppose that $G_n = (V_n, E_n)$ is either the hypercube $\Z_2^n$ or the
Caley graph of $S_n$ generated by transpositions. Then
\[
\var(|\CA_i|) = \tfrac{1}{4}\bigl(1+o(1)\bigr) |V_n|.
\]
\end{corollary}

In particular, the first-order asymptotics of the variance are exactly
the same as for an i.i.d. marking.

Throughout the remainder of the article, all graphs under consideration
shall satisfy Assumption~\ref{assumpmain}. In most examples, it will
be that $T_\mix^2(G_n) = o(F_{n,c})$ so that the second term in
(\ref{eqnvariance}) is negligible. In this case, taking $c = 2$ in
(\ref{eqnvariance}) provides a means to compute not only the magnitude but
also the constant in the first order asymptotics of the variance. In
some cases, such as $G_n = \Z_n^3$, the constant can even be computed
when $F_{n,c} = \Theta( (T_\mix(G_n))^2)$.

The challenge in obtaining Theorems~\ref{thmtori} and~\ref{thmmain}
is that the cancellation in the expansion of the variance is quite
significant which, when ignored, yields only an upper bound that can be
off by as much as a multiple of $T_\mix(G_n)$. We will now illustrate
this point in the case of $\Z_n^d$ for $d \geq3$. It will turn out
that the contribution to the variance from the sites visited by both
$X_1,X_2$ simultaneously is negligible, and hence we will ignore this
possibility in the present discussion. Observe
\begin{eqnarray*}
&&\var\biggl( \sum_{x} \one_{\{\tau_1(x) < \tau_2(x)\}} \biggr)\\
&&\qquad= \sum_{x,y}\bigl(\p[ \tau_1(x) < \tau_2(x), \tau_1(y) < \tau
_2(y)] \\
&&\qquad\hspace*{27pt}{}
-\p[\tau_1(x) < \tau_2(x)] \p[\tau_1(y) < \tau_2(y)]\bigr).
\end{eqnarray*}
Note that $\p[\tau_1(x) < \tau_2(x)]$ is approximately $\frac{1}{2}$.
Let $H(x,y) = \{ \tau_1(x) < \tau_1(y) \wedge\tau_2(x) \wedge\tau
_2(y)\}$. Consequently, by symmetry, the above is approximately equal to
\[
\sum_{x,y} \biggl(2 \p[ \tau_1(x) < \tau_2(x), \tau_1(y) < \tau_2(y)|
H(x,y)] \p[H(x,y)] - \frac{1}{4}\biggr) + O(n^d).
\]
The reason for the $O(n^d)$ term is that $\p[H(x,x)] = 0$, so all of
the diagonal terms are ignored in the summation.
Let $\wt{\pi}(\cdot;x,y)$ be the law of $X_2(\tau_1(x))$
conditional on
$H(x,y)$. As $\p[H(x,y)]$ is approximately $\frac{1}{4}$, using the
Markov property of $(X_1,X_2)$ applied for the stopping time $\tau
_1(x)$, we can rewrite the summation as
\[
2\sum_{x,y,z} \biggl(\p_{x,z}[\tau_1(y) < \tau_2(y)] - \frac{1}{2}
\biggr)\wt{\pi}(z;x,y) \p[H(x,y)].
\]
Here, $\p_{x,z}$ denotes the joint law of $X_1,X_2$ with $X_1(0) = x$
and $X_2(0) = z$. Thus we need to estimate
%
%
\begin{equation}
\label{eqnh14}
\frac{1}{2} \sum_{x,y,z} \biggl(\p_{x,z}[ \tau_1(y) <
\tau_2(y)] - \frac{1}{2} \biggr) \wt{\pi}(z;x,y).
\end{equation}
At this point, one is tempted to insert absolute values and then work
on each of the summands separately. Since $X_1$ and $X_2$ are
independent, note that $X_2(\tau_1(x)) \sim\pi(\cdot;\Z_n^d)$.
Thus by
Bayes' rule, we have
\[
\wt{\pi}(z;x,y) = \frac{\p[H(x,y)|X_2(\tau_1(x)) = z]}{\p
[H(x,y)]} \pi
(z;\Z_n^d) \leq C_0 \pi(z;\Z_n^d);
\]
see Theorem~\ref{thmrnestimate} for a much finer estimate. Hence the
expression in (\ref{eqnh14}) is bounded from above by
%
%
\begin{equation}
\label{eqnsumbound}
C_1 \sum_{x,y} \biggl| \p_{x,\pi}[\tau_1(y) < \tau_2(y)] - \frac{1}{2}
\biggr|,
\end{equation}
where $\p_{x,\pi}$ denotes the law of $X_1,X_2$ with $X_1(0) = x$ and
$X_2(0) \sim\pi(\cdot;\Z_n^d)$.

It is a basic fact that $T_\mix(\Z_n^d) = \Theta(n^2)$; one way to see
this is to invoke the local central limit theorem
(\cite{LAW91}, Theorem 1.2.1). We can analyze $\p_{x,\pi}[\tau_1(y)
< \tau
_2(y)]$ as follows. We consider two different cases: either $y$ is hit
before time $t_c \equiv cT_\mix(\Z_n^d) = c' n^2$ or afterward. The
probability that $X_2$ hits $y$ before $t_c$ is of order $n^{2-d}$ by a
union bound since $X_2(t) \sim\pi(\cdot;\Z_n^d) = n^{-d}$ for all $t$.
Second, by the local transience of random walk on $\Z_n^d$ for $d \geq
3$, the probability that $X_1$ hits $y$ before $t_c$ is, up to a
multiplicative constant, well approximated by $g(x,y;\Z_n^d)$. We now
consider the second case. By time $t_c$ for $c > 0$ large enough, $X_1$
will have mixed. This means that if neither $X_1$ nor $X_2$ has hit $y$
by this time, the probability that either one hits first is close to
$1/2$. The careful reader who wishes to see precise, quantitative
versions of these statements will find such in the lemmas we use to
prove Theorem~\ref{thmmain}. Thus it is not difficult to see that
there exists $C_2 > 0$ so that
\[
|\p_{x,\pi}[\tau_1(y) < \tau_2(y)] - 1/2| \leq C_2 g(x,y;\Z_n^d).
\]
This leads to an upper bound of
\[
C_3 \sum_{x,y} g(x,y;\Z_n^d) \leq C_4 n^{d+2}.
\]
A slightly more refined analysis leads to a lower bound of (\ref
{eqnsumbound}) with the same growth rate. As we will show in the next
section, in every dimension this estimate is typically quite far from
being sharp. The reason for the inaccuracy is that by moving the
absolute value into the sum in (\ref{eqnsumbound}) we are unable to
take advantage of the cancellation that arises as $\p_{x,\pi}[ \tau
_1(y) < \tau_2(y)] > 1/2$ when $x$ is close to $y$ and $\p_{x,\pi}[
\tau
_1(y) < \tau_2(y)] < 1/2$ when $x$ is far from $y$.

\subsection*{Outline} The remainder of this article is structured as
follows. In the next section, we will deduce Theorem~\ref{thmtori}
and Corollary~\ref{corhypsym} from Theorem~\ref{thmmain}. In
Section~\ref{secprelim}, we introduce some notation that will be used
throughout in addition to collecting several basic random walk
estimates. Next, in Section~\ref{secrn}, we give a precise estimate
of the Radon--Nikodym derivative of $\wt{\pi}(\cdot;x,y)$ with respect
to $\pi$. In Section~\ref{secvar}, we prove Theorem~\ref{thmmain}
and end in Section~\ref{secproblems} with a list of related problems
and discussion.

\section{\texorpdfstring{Proof of Theorem \protect\ref{thmtori} and Corollary \protect\ref{corhypsym}}
{Proof of Theorem 1.1 and Corollary 1.4}}

The following proposition will be important for the proof of Theorem
\ref{thmtori}.
%
%
\begin{proposition}
\label{proplimitexistence}
Assume that $G_n = \Z_n^d$ for $d \geq3$. For each $c > 1$, the limit
%
%
\begin{equation}
\label{eqnlimitexistence}
\lim_{n \to\infty} \frac{1}{h_d(n)} \sum_{x,y} \bigl(f_{n,c}(x,y) - \ol
{f}_{n,c}\bigr)^2
\end{equation}
exists. When $d \geq4$, it is $\alpha_d$ as in (\ref{eqnalphad5}),
(\ref{eqnalphad4}). When $d=3$, it is given by $\alpha_3^{c}$ where
$\alpha_3^T$ is as in (\ref{eqnalphad3}).
\end{proposition}

The first step in the proof of the proposition is to reduce the
existence of the limit to a computation involving Green's functions.
Recall from (\ref{eqngreensfunctionzd}) that $G(y;\Z^d)$ is the
Green's function for lazy random walk on $\Z^d$. In order to keep the
notation from becoming too heavy, throughout the rest of this section
we will write $T_\mix$ for $T_\mix(G_n)$ where $G_n$ will be clear from
the context. Let
\[
g_c(x,y;G_n) = \E_x \sum_{t=0}^{c T_\mix} \one_{\{X(t) = y\}}.
\]

%
\begin{lemma}
\label{lemgreenreduction}
Assume that $G_n = \Z_n^d$ for $d \geq3$. For each $c > 1$, we have that
\[
\lim_{n \to\infty} \frac{1}{h_d(n)} \sum_{x,y} \bigl(f_{n,c}(x,y) - G^{-1}
(0;\Z^d) g_c(x,y;\Z_n^d)\bigr)^2 = 0.
\]
\end{lemma}
\begin{pf}
Observe
\[
g_c(x,y;\Z_n^d) \leq f_{n,c}(x,y) g_c(y,y;\Z_n^d).\vadjust{\goodbreak}
\]
We shall now prove a matching lower bound. Fix $0 < \wt{c} < c$. Then
we have that
%
%
\begin{eqnarray}\label{eqngreenlowerbound}
g_c(x,y;\Z_n^d)
&\geq&\E_x \Biggl[ \Biggl(\sum_{t=\tau(y)}^{c T_\mix} \one_{\{X(t)
= y\}
}\Biggr) \one_{\{\tau(y) \leq(c-\wt{c}) T_\mix\}} \Biggr]
\nonumber\\[-8pt]\\[-8pt]
&\geq& f_{n,c-\wt{c}}(x,y)g_{\wt{c}}(y,y;\Z_n^d).\nonumber
\end{eqnarray}
Assuming $c - \wt{c} > 1$, by mixing considerations as well as a union
bound (see Proposition~\ref{propmixingdecay}) we have that
%
%
\begin{eqnarray} \label{eqnhitmixbound}
f_{n,c-\wt{c}}(x,y)
&=& f_{n,c}(x,y) - \p_x[(c-\wt{c}) T_\mix< \tau(y) \leq cT_\mix]
\nonumber\\[-8pt]\\[-8pt]
&=& f_{n,c}(x,y) + O(\wt{c} n^{2-d}).\nonumber
\end{eqnarray}
Since $\wt{c} > 0$, we have
%
%
\begin{eqnarray}
\label{eqngreendecay}
g_{\wt{c}}(y,y;\Z_n^d) &=& g_c(y,y;\Z_n^d) - \sum_{z} p^{\wt{c}
T_\mix
}(y,z;\Z_n^d) g_{c-\wt{c}}(z,y;\Z_n^d)\nonumber\\[-8pt]\\[-8pt]
&=& g_c(y,y;\Z_n^d) + O\bigl((c-\wt{c})\wt{c}^{-d/2} n^{2-d}\bigr), \nonumber
\end{eqnarray}
where we used in the last line that $p^t(z,y;\Z_n^d) \leq c_1 t^{-d/2}$
for some $c_1 > 0$ (see~\cite{LAW91}, Theorem 1.2.1) as well as the
observation $\sum_z g_{c-\wt{c}}(z,y;\Z_n^d) = (c-\wt{c}) T_\mix$.
Combining (\ref{eqngreenlowerbound}), (\ref{eqnhitmixbound})
and (\ref{eqngreendecay}), we have thus proved the lower bound
\[
g_c(x,y;\Z_n^d) \geq f_{n,c}(x,y) g_c(y,y;\Z_n^d) + O\bigl((c-\wt{c})\wt
{c}^{-d/2}\bigl(\wt{c} + g_c(x,y;\Z_n^d)\bigr) n^{2-d}\bigr).
\]
Here, we used the bound $f_{n,c}(x,y) \leq g_c(x,y;\Z_n^d)$.
Theorem 1.5.4 of~\cite{LAW91} implies $g_c(x,y;\Z_n^d) = \Theta(c
|x-y|^{2-d})$ (it is actually stated for walks on $\Z^d$ which are not
lazy, but the generalization is straightforward). Consequently,
\[
\sum_{y} g_c^2(x,y;\Z_n^d) = \cases{
\Theta(n), &\quad if $d = 3$,\cr
\Theta(\log n), &\quad if $d = 4$, and\cr
\Theta(1), &\quad if $d \geq5$.}
\]
Hence,
\begin{eqnarray*}
&&\sum_{x,y} \bigl(f_{n,c}(x,y) g_c(y,y;\Z_n^d) - g_c(x,y;\Z_n^d)\bigr)^2\\
&&\qquad= \sum_{x,y} \bigl[O\bigl((c-\wt{c}) \wt{c}^{-d/2}\bigl(\wt{c} + g(x,y;\Z_n^d)\bigr)
n^{2-d}\bigr)\bigr]^2\\
&&\qquad= O\bigl((c-\wt{c})^2 \wt{c}^{-d} \bigl(\wt{c}^2 + o(1)\bigr) n^4\bigr).
\end{eqnarray*}
Dividing both sides by $h_d(n)$, taking a limsup as $n \to\infty$,
then as $\wt{c} \to0$ yields
\[
\lim_{n \to\infty} \frac{1}{h_d(n)} \sum_{x,y} \bigl(f_{n,c}(x,y)
g_c(y,y;\Z
_n^d) - g_c(x,y;\Z_n^d)\bigr)^2 = 0.
\]
By (\ref{eqngreendecay}) we\vspace*{1pt} know that $|g_c(y,y;\Z_n^d) -
g_1(y,y;\Z
_n^d)| = o(1)$, and, by local transience, it is not hard to see that
$\lim_{n \to\infty} g_1(y,y;\Z_n^d) = G(0;\Z^d)$.
\end{pf}
\begin{pf*}{Proof of Proposition~\ref{proplimitexistence}}
Lemma~\ref{lemgreenreduction} implies that we may\vspace*{1pt} replace
$f_{n,c}(x,y)$ by $G^{-1}(0;\Z^d) g_c(x,y;\Z_n^d)$ in (\ref
{eqnlimitexistence}). Letting $\ol{g}_{n,c} = c T_\mix n^{-d}$, we
can likewise replace $\ol{f}_{n,c}$ in (\ref{eqnlimitexistence}) by
$G^{-1}(0;\Z^d) \ol{g}_{n,c}$. Consequently, to prove the proposition,
it suffices to prove the existence of the limit
%
%
\begin{equation} \label{eqngreenlimit}
\lim_{n \to\infty} \frac{1}{h_d(n)} \sum_{x,y} \bigl(g_c(x,y;\Z_n^d) -
\ol
{g}_{n,c}\bigr)^2.
\end{equation}
We will divide the proof into the cases $d \geq4$ and $d=3$.

\textit{Case} 1: $d \geq4$. As $\ol{g}_{n,c} = O(c n^{2-d})$,
we have
\[
\frac{1}{h_d(n)} \sum_{x,y} \bigl(\ol{g}_{n,c}^2 + 2 \ol{g}_{n,c}
g_c(x,y;\Z
_n^d)\bigr) = o(1).
\]
Thus it suffices to show in this case that
\[
\lim_{n \to\infty} \frac{1}{\wt{h}_d(n)} \sum_{y} g_c^2(0,y;\Z_n^d)
\]
exists, where $\wt{h}_4(n) = \log n$ and $\wt{h}_d(n) = 1$ for $d
\geq
5$. This will be a consequence of two observations. First, note that
\begin{eqnarray*}
\sum_{|y| \geq\ell} g_c^2(0,y;\Z_n^d)
&=& \sum_{|y| \geq\ell} O(c |y|^{4-2d})
= \sum_{m=\ell}^n O(c m^{4-2d} \cdot m^{d-1})\\
&=& \cases{
O(c \ell^{-1}), &\quad if $d \geq5$,\cr
O\bigl(c \log(n/ \ell)\bigr), &\quad if $d = 4$.}
\end{eqnarray*}
Thus it suffices to show that, for $\ell=
\ell(n,\varepsilon) = n^{1-\varepsilon}$ with $\varepsilon> 0$, the limit
\[
\lim_{\varepsilon\to0} \lim_{n \to\infty} \frac{1}{\wt{h}_d(n)}
\sum
_{|y| \leq\ell} g_c^2(0,y;\Z_n^d)
\]
exists (we can even restrict to finite $\ell$ if $d \geq5$). Our
second observation is that
\[
g_c(0,y;\Z_n^d) - G(y;\Z^d) = O(c n^{2-d}) \qquad\mbox{for } |y| \leq
\ell.
\]
This follows since we can couple the walks on $\Z_n^d$ and $\Z^d$
starting at $0$ such that they are the same until the first time $\tau
_0$ they have reached distance $n/2$ from $0$, then move independently
thereafter. The expected number of visits each walk makes to $y$ after
time $\tau_0$, where the former is stopped at time $c T_\mix$, is
easily seen to be $O(c n^{2-d})$.
Thus,
\[
\sum_{|y| \leq\ell} \bigl(g_c(0,y;\Z_n^d) - G(y;\Z^d)\bigr)^2 = o(1).
\]
Therefore if $d \geq5$, we have
\[
\lim_{n \to\infty} \frac{1}{h_d(n)}\sum_{x,y} g_c^2(x,y;\Z_n^d) =
\sum
_{y \in\Z^d} G^2(y;\Z^d).
\]
For $d = 4$,
\[
\lim_{n \to\infty} \frac{1}{h_4(n)}\sum_{x,y} g_c^2(x,y;\Z_n^4) =
\lim
_{n \to\infty} \frac{1}{\log n}\sum_{y \in\Z^4, |y| \leq n}
G^2(y;\Z^4).
\]
Note that the limit on the right-hand side exists since by Theorem
1.5.4 of~\cite{LAW91} (generalized to lazy walks)
\[
G(y;\Z^d) = a_d |y|^{2-d} + o(|y|^{-\alpha}),
\]
if $\alpha\in(0,d)$ is fixed.

\textit{Case} 2: $d=3$. The thrust of the previous argument was
that random walk on $\Z_n^d$ for $d \geq4$ is sufficiently transient
so that pairs of points of distance $\Omega(n^{1-\varepsilon})$ make a
negligible contribution to the variance, which in turn allowed us to
make an accurate comparison between the Green's function for random
walk on $\Z_n^d$ with that on $\Z^d$. The situation for $d=3$ is more
delicate since the opposite is true: pairs of distance $O(n^{1-\varepsilon
})$ do not measurably affect the variance.

Theorem 1.2.1 of~\cite{LAW91} (extended to the case of lazy random
walk, see also Corollary 22.3 of~\cite{BR76}) implies the existence of
constants $\beta_3,\gamma_3 > 0$ such that with $\ol{p}^t(x,y;\Z^3) =
\frac{\beta_3}{t^{3/2}} \exp( - \frac{\gamma_3|x-y|^2}{t}
)$, we have the estimate
\[
|\ol{p}^t(x,y;\Z^3) - p^t(x,y;\Z^3)| = |x-y|^{-2}O(t^{-3/2}).
\]
Hence letting $\ol{p}^t(x,y;\Z_n^3) = \sum_{k \in\Z^3} \ol
{p}^t(x,y+kn;\Z^3)$, one can easily show that with
\[
\Delta(x,y) \equiv\sum_{t=0}^{c T_\mix} |\ol{p}^t(x,y;\Z_n^3) -
p^t(x,y;\Z_n^3)|
\]
we have that
%
%
\begin{equation}
\label{eqngreenboundsum}
\frac{1}{h_3(n)} \sum_{x,y} \Delta^2(x,y) = o(1).
\end{equation}
By differentiating $\ol{p}$ in $t$, we see that for $1 \leq t \leq s
\leq t+1$, we have
\begin{eqnarray*}
&&|\ol{p}^s(0,y;\Z_n^3) - \ol{p}^t(0,y;\Z_n^3)| \\
&&\qquad=
O\biggl(\frac{\ol{p}^t(0,y;\Z_n^3)}{t} + \sum_k \frac{|y+kn|^{2}}{t^2}
\ol{p}^t(0,y+kn;\Z_n^3)\biggr) .
\end{eqnarray*}
We are now going to prove that
%
%
\begin{equation}
\label{eqngreenboundtime}
\sum_{y \in\Z_n^3} \biggl( \int_1^{c
T_\mix} \ol{p}^t(0,y;\Z_n^3) - \ol{p}^{\lfloor t
\rfloor}(0,y;\Z_n^3) \,dt \biggr)^2 = O(1).
\end{equation}
It suffices to bound
\begin{eqnarray*}
A &\equiv& \sum_{y \in\Z_n^3} \biggl( \int_1^{c T_\mix} \frac
{1}{t} \ol
{p}^t(0,y;\Z_n^3) \,dt \biggr)^2,\\
B &\equiv& \sum_{y \in\Z_n^3} \biggl(\sum_k \int_1^{c T_\mix}
\frac{|y+kn|^{2}}{t^2}
\ol{p}^t(0,y+kn;\Z_n^3)
\,dt \biggr)^2.
\end{eqnarray*}
For $A$, we apply Cauchy--Schwarz to the integral and invoke the
integrability of $1/t^2$ over $[1,\infty)$ to arrive at
\[
A \leq C_2 \sum_{y \in\Z_n^3} \int_1^{c T_\mix}
[\ol{p}^t(0,y;\Z_n^3)]^2 \,dt = O(1).
\]
For $B$, we insert the formula for $\ol{p}$ into the integral, make the
substitution $u=|y+kn|^{2}/t$ and then compute to see
\[
B \leq C_3 \sum_{y \in\Z_n^3} \frac{1}{|y|^6+1} = O(1).
\]
This proves (\ref{eqngreenboundtime}). Recall that $\T^3$ is the
three-dimensional continuum torus. For $x,y \in\T^3$, let
%
%
\begin{equation}
\label{eqngreenfunctionbmtorus}
g_c(x,y;\T^3)
= \int_{0}^{c T_\mix} \ol{p}^t(nx,ny;\Z_n^3) \,dt = \frac{1}{n}
\int_0^{cT} \ol{p}^{u}(x,y;\T^3)\,du,
\end{equation}
where $T = T_\mix/ n^2$.
By (\ref{eqngreenboundsum}), (\ref{eqngreenboundtime}), we
have that
\[
\frac{1}{h_3(n)} \sum_{x,y \in\Z_n^3} \bigl(g_c(x,y;\Z_n^3) -
g_c(x/n,y/n;\T
^3)\bigr)^2 = o(1).
\]
Therefore we may replace $g_c(x,y;\Z_n^3)$ in (\ref{eqngreenlimit})
with $g_c(x/n,y/n;\T^3)$. Note that $g_c(\cdot,\cdot;\T^3)$ is the
product of $n^{-1}$, and the Green's function for $B_{t/2}$, where $B_t$
is a Brownian motion on $\T^3$; roughly, the reason that the Brownian
motion moves at $1/2$-speed is that a lazy random walk moves at $1/2$
the speed of a simple random walk. It is left to bound
\[
n^2 \int_{\T^3} \int_{\T^3} \bigl( g_c( \lfloor nx \rfloor/ n, \lfloor ny
\rfloor/ n;\T^3) - g_c(x,y;\T^3)\bigr)^2 \,dx \,dy;
\]
the reason for the pre-factor $n^2$ is that we need to multiply by
$(n^3)^2$ in order to make the double integral comparable to the double
summation, and we also divide by the normalization $h_3(n)$. From
(\ref{eqngreenfunctionbmtorus}), we see that $g_c(x,y;\T^3)$ is
$O(n^{-1})$-Lipschitz away from the diagonal $D_\varepsilon= \{ (x,y) \in
\T^3 \times\T^3\dvtx|x-y| \leq\varepsilon\}$. Thus since $|(x,y) - (
\lfloor nx \rfloor/n , \lfloor ny \rfloor/n)| = O(n^{-1})$, the
integrand is $O(n^{-4})$ on $D_\varepsilon^c$, hence the integral\vspace*{1pt} over
$D_\varepsilon^c$ is $O(n^{-2})$. Since both $n g_c( \lfloor nx \rfloor/
n, \lfloor ny \rfloor/ n;\T^3)$ and $n g_c(x,y;\T^3)$ are uniformly
$L^2$-integrable over $\T^3 \times\T^3$, it follows that the
contribution coming from $D_\varepsilon$ can be made uniformly small in
$n$ by first fixing $\varepsilon> 0$ small enough.
\end{pf*}

We now deduce Theorem~\ref{thmtori} from Theorem~\ref{thmmain}.
\begin{pf*}{Proof of Theorem~\ref{thmtori}}
Suppose $G_n = \Z_n^d$ for $d \geq3$. Recall that
\[
T_\mix(\Z_n^d) =
\Theta(n^2)
\]
(see~\cite{LPW08}) and there exists $c_d > 0$ so that
$g(x,y;\Z_n^d) \leq c_d|x-y|^{2-d} \wedge1$ (see~\cite{LAW91}).
Consequently, the hypotheses of Theorem~\ref{thmmain} are obviously
satisfied, except for possibly~\ref{assumpmainhitmiss}. This is
easy to see if $x$ is sufficiently far from $y,z$ so that
$g(x,y;\Z_n^d) + g(x,z;\Z_n^d) \leq1/2$. Now suppose that $|x-y|
\wedge|x-z| = r$ is small enough so that $g(x,y;\Z_n^d) +
g(x,z;\Z_n^d) > 1/2$. We have the trivial bound that $X$ starting at
$x$ will get to distance $r+s$ without hitting $y,z$ in $s$ steps with
probability at least $(4d)^{-s}$ since in each step, $X$ has
probability at least $(4d)^{-1}$ of increasing its distance from $y,z$
by $1$. If $s$ is large enough, then after such steps we will have
$g(X_s,y;\Z_n^d) + g(X_s,z;\Z_n^d) \leq1/2$, which gives the desired
result.

Proposition~\ref{proplimitexistence} implies that $F_{n,c} \sim
\frac{1}{4}\alpha_{d,c} h_d(n)$ as $n \to\infty$. This is enough to
dominate $T_\mix^2(\Z_n^d) = \Theta(n^4)$ except if $d=3$. We shall now
argue that, nevertheless, $F_{n,c}$ is still the dominant term in this case.
Note that
\[
\ol{f}_{n,c} \leq\frac{1}{n^3} \sum_y g_c(x,y;\Z_n^3) \leq A_0 c n^{-1}
\]
for some $A_0 > 0$ and $c \geq2$ fixed. Also, the transience of random
walk on $\Z_n^3$ implies that there exists $A_1 > 0$ so that
$f_{n,c}(x,y) \geq A_1 |x-y|^{-1} \wedge1$. Thus for
\[
|x-y| \leq\biggl(\frac{A_1}{2 A_0 c}\biggr) n \equiv A_2 n
\]
we have that $f_{n,c}(x,y) - \ol{f}_{n,c} \geq\frac{A_1}{2}|x-y|^{-1}
\wedge1$. Consequently,
\[
F_{n,c} \geq\frac{A_1^2}{4} \sum_{|x-y| \leq A_2 n} |x-y|^{-2}
\wedge
1 = c^{-1} \Theta(n^4).
\]
A matching upper bound, up to a multiplicative factor, is also not
difficult to see.

Our lower bound for $F_{n,c}$ depends on $c$ by a multiplicative factor
of $1/c$ while the second term in (\ref{eqnvariance}) decays
exponentially in $c$. Thus by taking $c \geq2$ large enough, we see
that $F_{n,c}$ is still dominant for $d = 3$.
\end{pf*}

We now turn to the proof of Corollary~\ref{corhypsym}.
\begin{pf*}{Proof of Corollary~\ref{corhypsym} for the Hypercube}
For $\Z_2^n$, it is easier to work with the continuous time random
walk (CTRW) since the types of estimates we require easily translate
over to the corresponding lazy walk. The transition kernel of the CTRW is
\[
p^t(x,y;\Z_2^n) = \frac{1}{2^n} (1+e^{-2t/n})^{n-|x-y|} (1-e^{-2t/n})^{|x-y|},
\]
where $|x-y|$ is the number of coordinates in which $x$ and $y$ differ.
The spectral gap is $1/n$ (see Example 12.15 of~\cite{LPW08}) which
implies $\Omega(n) = T_\mix(\Z_2^n) = O(n^2)$ (see Theorem 12.3 of
\cite{LPW08}). Consequently, the first hypothesis of Theorem~\ref{thmmain}
holds. If $|x-y| = r$, then it is easy to see there exists $C_\varepsilon
,\rho_\varepsilon> 0$ so that
\[
p^t(x,y;\Z_2^n) \leq\cases{
\displaystyle \biggl(C_\varepsilon\frac{t}{n}\biggr)^r \exp\biggl(- \frac{t}{C_\varepsilon
n}(n-r)\biggr),
&\quad if $t \leq\varepsilon n$,\cr
e^{-\rho_\varepsilon n}, &\quad if $t > \varepsilon n$,}
\]
provided $\varepsilon> 0$ is sufficiently small. Thus it is not difficult
to see that $g(x,y;\Z_2^n) \leq C_\varepsilon' n^{-r}$.
Trivially,
\[
\bigl|\{y \in\Z_2^n\dvtx|x-y| = r\}\bigr| = \pmatrix{ n \cr r} \leq n^r.
\]
Thus for $x_0$ fixed we have
\[
\sum_{y \neq x_0} g^2(x_0,y;\Z_2^n) \leq O\Biggl( \sum_{r=1}^n n^{-2r}
\cdot n^r \Biggr) = O\biggl( \frac{1}{n} \biggr),
\]
so the second hypothesis of Theorem~\ref{thmmain} is satisfied. The
final hypothesis is obviously also satisfied. Now, a union bound
implies that $\ol{f}_{n,c} = O(2^{-n} T_\mix(\Z_2^n))$, which implies
$(f_{n,c}(x,x) - \ol{f}_{n,c})^2 = 1+o(1)$. On the other hand,
\[
\sum_{|x-y| \geq1} f_{n,c}^2(x,y) = O\Biggl( 2^n \sum_{r=1}^n n^{-2r}
\cdot n^r \Biggr) = o(2^n).
\]
Putting everything together, Theorem~\ref{thmmain} implies
\[
\var( |\CA_i| ) = \tfrac{1}{4}\bigl(1+o(1)\bigr) 2^n.
\]
\upqed\end{pf*}
\begin{pf*}{Proof of Corollary~\ref{corhypsym} for the Caley graph of
$S_n$}
Let $G_n$ be the Caley graph of $S_n$ generated by transpositions. By
work of Diaconis and Shashahani~\cite{DS81}, the total variation mixing
time of $S_n$ is $\Theta(n \log n)$, which by Theorem~12.3 of
\cite{LPW08} implies $T_\mix(G_n) = O( n (\log n) (\log n!)) = O(n^2
(\log n)^2)$. We are now going to give a crude estimate of
$p^t(\sigma,\tau;S_n)$. By applying an automorphism, we may assume
without loss of generality that $\sigma= \mathrm{id}$. Suppose that
$d(\mathrm{id},\tau) = r$ and that $\tau_1,\ldots,\tau_r$ are
transpositions such that $\tau_r \cdots\tau_1 = \tau$. Then
$\tau_1,\ldots,\tau_r$ move at most $2r$ of the $n$ elements of
$\{1,\ldots,n\}$, say, $k_1,\ldots,k_{2r}$. Suppose
$k_1',\ldots,k_{2r}'$ are distinct from $k_1,\ldots,k_{2r}$ and
$\alpha
\in S_n$ is such that $\alpha(k_i) = k_i'$ for $1 \leq i \leq r$. Then
the automorphism of $G_n$ induced by conjugation by $\alpha$ satisfies
$\alpha\tau\alpha^{-1} \neq\tau$. Therefore the size of the set of
elements $\tau'$ in $S_n$ such that there exists a graph automorphism
$\varphi$ of $G_n$ satisfying $\varphi(\tau) = \tau'$ and $\varphi
(\mathrm{id}) = \mathrm{id}$ is at least ${n\choose2r} \geq
2^{-2r} n^{2r} ((2r)!)^{-1}$, assuming $n \geq4r$. Therefore,
%
%
\begin{equation}
\label{eqnpermbound}
p^t(e,\tau;G_n) \leq\frac{2^{2r} (2r)!}{n^{2r}}
\quad\mbox{and}\quad
g(e,\tau
;G_n) \leq C(2^{2r} (2r)!) (\log n)^2 n^{2-2r}.\hspace*{-26pt}
\end{equation}
This bound is good enough for $r \geq2$, but does not quite suffice
when $r = 1$. This case is not difficult to handle, however, since it
is easy to see that the random walk has distance $3$ from $e$ with
probability $1 - O(1/n)$ after its first three moves, hence with
distance at least $2$ from any permutation with distance $1$ from $e$.
Combining this with (\ref{eqnpermbound})
gives a bound on $g(e,\tau ;G_n)$ for all $\tau\in S_n$.
From this is it clear that $(G_n)$ satisfies the hypotheses of Theorem
\ref{thmmain} and, arguing as in the case of the hypercube, that
\[
\var( |\CA_i| ) = \tfrac{1}{4}\bigl(1+o(1)\bigr) n!.
\]
\upqed\end{pf*}

\section{Preliminaries}
\label{secprelim}

\subsection{Notation}
Suppose that $G = (V,E)$ is a graph, and let $X_1,X_2$ be independent
random walks on $G$. Recall that $a \wedge b = \min(a,b)$ for $a,b \in
\R$. For $x,y \in V$, let
\[
\tau_i(x,y) = \tau_i(x) \wedge\tau_i(y) \quad
\mbox{and}\quad \tau(x,y) =
\tau
_1(x,y) \wedge\tau_2(x,y),
\]
where $\tau_i(x) = \min\{ t \geq0\dvtx X_i(t) = x\}$.
Let
\[
H(x,y) =
\{ \tau_1(x) < \tau_1(y) \wedge\tau_2(x,y)\}.
\]
This is the event that $x$ is hit by $X_1$ before $X_2$ as well as
before both $X_1,X_2$ hit $y$. Let
\[
\wt{\pi}(z;x,y) = \p[ X_2(\tau_1(x,y)) = z | H(x,y)],
\]
and let $\pi$ be the uniform measure on $V$. Throughout, $\p_z[\cdot]$
denotes the law of random walk initialized at $z$ (and the initial
distribution is stationary whenever $z$ is omitted).\vadjust{\goodbreak}
The proofs in this article will involve probabilities of complicated
events. To keep the formulas succinct, it will be helpful for us to
introduce the following notation: let
\begin{eqnarray*}
G_{ij}(x) &=& \{ \tau_i(x) < \tau_j(x)\},\\
G_{ij}(x,y) &=& \{ \tau_i(x,y) < \tau_j(x,y)\}
\end{eqnarray*}
and
\[
G_{i}(x,y) = \{\tau_i(x) < \tau_i(y)\}.
\]
Throughout we will fix a sequence of graphs $(G_n)$ satisfying
Assumption~\ref{assumpmain}. We let
\begin{eqnarray*}
\Gamma_n &=& {c_0 T_\mix(G_n) \log}|V_n|,\qquad
\Upsilon_n = \frac{T_\mix(G_n)}{ |V_n|},\\
\Delta_n &=& {\Upsilon_n \log}|V_n|,\qquad
\CS_n = \sum_{y \neq x_0} g^2(x_0,y;G_n),
\end{eqnarray*}
where $c_0$ will be determined later, and $x_0$ is fixed. Note that
$\CS
_n$ does not depend on $x_0$ by vertex transitivity. We will typically
write $T_\mix$ for $T_\mix(G_n)$, $p^t(\cdot,\cdot)$ for $p^t(\cdot
,\cdot;G_n)$ and $g(\cdot,\cdot)$ for $g(\cdot,\cdot;G_n)$ in
order to
keep the notation light and, in general, suppress dependencies on $n$.

\subsection{Elementary estimates}
Recall that the total variation distance of probability measures $\mu
,\nu$ on $V$ is
\[
\| \mu-\nu\|_{\mathrm{TV}} = {\max_{A \subseteq V}} | \mu(A) - \nu(A)| =
\frac
{1}{2} \sum_{x \in V} |\mu(\{x\}) - \nu(\{x\})|.
\]
The following provides a bound on the rate of decay of the distance of
$p^t(x,\cdot)$ to stationarity.
%
%
\begin{proposition}
\label{propmixingdecay}
For every $s,t \in\N$,
%
%
\begin{eqnarray}\quad
\label{eqntvdecay}
{\max_x }\| p^{t+s}(x,\cdot) - \pi\|_{\mathrm{TV}} &\leq& {4 \max_{x,y} }\|
p^t(x,\cdot) - \pi\|_{\mathrm{TV}} \| p^s(y,\cdot) -
\pi\|_{\mathrm{TV}},
\\
\label{eqnuniformdecay}
{\max_{x,y}} \biggl| \frac{p^{t+s}(x,y)}{\pi(y)} - 1\biggr| &\leq&
\max_{x,y} {\frac{p^{s}(x,y)}{\pi(y)} \max_x} \| p^{t}(x,\cdot) - \pi\|_{\mathrm{TV}}.
\end{eqnarray}
\end{proposition}
\begin{pf}
The first part is a standard result; see, for example, Lemmas~4.11
and~4.12 of~\cite{LPW08}. The second part is a consequence of the
semigroup property:
\begin{eqnarray*}
\frac{1}{\pi(z)} p^{t+s}(x,z)
&=& \frac{1}{\pi(z)} \sum_y p^t(x,y) p^s(y,z)\\
&=& \frac{1}{\pi(z)} \sum_y [p^t(x,y) - \pi(y) + \pi(y)]p^s(y,z)\\
&\leq& \biggl(\max_{y,z} \frac{p^s(y,z)}{\pi(z)}\biggr) \|
p^t(x,\cdot)
- \pi\|_{\mathrm{TV}} + 1.\hspace*{15pt}\qed
\end{eqnarray*}
\noqed\end{pf}

Trivially,
\[
{\max_x }\| p^t(x,\cdot) - \pi\|_{\mathrm{TV}} \leq\max_{x,y} \biggl| \frac
{p^t(x,y)}{\pi(y)} - 1 \biggr|.
\]
Consequently, (\ref{eqntvdecay}) and (\ref{eqnuniformdecay}) give
%
%
\begin{equation}\label{eqnuniformexponentialdecay}
\max_{x,y} \biggl| \frac{p^{t}(x,y)}{\pi(y)} - 1 \biggr| \leq
\ol{\gamma} e^{-\ol{\gamma}\alpha} \qquad\mbox{for } t \geq\alpha
T_\mix
\mbox{ and } \alpha> 0,
\end{equation}
where $\ol{\gamma} > 0$ is a universal constant. We will often use
(\ref{eqnuniformexponentialdecay}) without reference.

Throughout the article, it will be important for us to have precise
estimates of the Radon--Nikodym derivative of the law of random walk
conditioned on various events with respect to the uniform measure. In
the following, we are interested in the case of a random walk
conditioned not to have hit a particular point. Let $T_k = k T_\mix$.
%
%
\begin{lemma}
\label{lemconditionnothitrn}
There exists $\gamma,p_0 > 0$ so that for all $k \geq1$ satisfying $k
\Upsilon_n \leq p_0$ and $c \geq2$, we have
\[
\p_x[X(c T_k) = z|\tau(y) > c T_k]
= \bigl[1 + O\bigl(e^{-\gamma c k } + c k \Upsilon_n + g(y,z)\bigr)\bigr] \pi(z).
\]
\end{lemma}

Note that by part~\ref{assumpmainmixingbound} of Assumption
\ref{assumpmain}, this lemma applies if
\[
k = O( ({\log}|V_n|)^2).
\]
\begin{pf*}{Proof of Lemma~\ref{lemconditionnothitrn}}
Using that $\p_x[ X(c T_k) = z] = (1+O(e^{-\ol{\gamma} ck})) \pi(z)$,
an application of Bayes' formula yields
\begin{eqnarray*}
&& \p_x[X(c T_k) = z|\tau(y) > c T_k]\\
&&\qquad= \frac{\p_x[ \tau(y) > c T_k | X(c T_k) = z]}{\p_x[ \tau(y) > c T_k]}
\bigl(1+O(e^{-\ol{\gamma} c k})\bigr) \pi(z).
\end{eqnarray*}
The idea of the rest of the proof is to show it is unlikely that $X$
hits $y$ close to time $cT_k$, in which case we can use a mixing
argument to show that conditioning on $X(c T_k) = z$ has little effect.
For $1 \leq\wt{c} \leq\wt{c}+1 \leq c$, we have
\begin{eqnarray*}
&&\p_x[ \tau(y) > c T_k | X(c T_k) = z]\\
&&\qquad=\p_x[ \tau(y) > \wt{c} T_k | X(c T_k) = z] - \p_x[ c T_k \geq
\tau
(y) > \wt{c} T_k | X(c T_k) = z].
\end{eqnarray*}
By a time-reversal, we have that
\[
\p_x[c T_k \geq\tau(y) > \wt{c} T_k | X(cT_k) = z] \leq\p_z[ \tau(y)
< (c-\wt{c}) T_k | X(cT_k) = x].
\]
By mixing considerations and a union bound, we have
\[
\p_z[ \tau(y) \leq(c-\wt{c}) T_k | X(c T_k) = x] =
O\bigl( g(y,z) + (c-\wt{c}) k \Upsilon_n\bigr).
\]
Applying Bayes' formula, observe
\begin{eqnarray*}
\p_x[ \tau(y) > \wt{c} T_k | X(c T_k) = z]
&=& \frac{\p_x[ X(c T_k) = z | \tau(y) > \wt{c} T_k]}{\p_x[X( cT_k) =
z]} \p_x[\tau(y) > \wt{c} T_k]\\
&=& \bigl(1+O\bigl(e^{-k \ol{\gamma}(c-\wt{c})}\bigr)\bigr)\p_x[\tau(y) > \wt{c} T_k].
\end{eqnarray*}
Similarly,
\[
\p_x[ \tau(y) > c T_k]
= \p_x[ \tau(y) > \wt{c} T_k] - \p_x[ c T_k \geq\tau(y) > \wt{c} T_k].
\]
By a union bound and mixing considerations, the second term on the
right-hand side is of order $O((c-\wt{c})k \Upsilon_n)$.
We now take $\wt{c} = c/2$ and $\gamma= \ol{\gamma}/2$. By part
\ref{assumpmainhitmiss} of Assumption~\ref{assumpmain}, we have that
\[
\p_x[\tau(y) > cT_k] \geq1 - \rho_0 - O(c k \Upsilon_n)
\]
uniformly in $n$. In particular, there exists $p_0 > 0$ so that if $c k
\Upsilon_n \leq p_0$, then $\p_x[\tau(y) > cT_k]$ is uniformly positive
in $n$. Putting everything together, for such~$k$, we thus have
\begin{eqnarray*}
&& \p_x[ X(c T_k) = z | \tau(y) > c T_k]\\
&&\qquad= \frac{(1+ O(e^{-\gamma c k}))\p_x[\tau(y) > \wt{c} T_k] +
O(g(y,z) +
c k \Upsilon_n)}{\p_x[\tau(y) > \wt{c} T_k] + O(c k \Upsilon_n)}\\
&&\qquad\quad{}\times\bigl(1+O(e^{-\ol{\gamma} c k})\bigr)\pi(z)\\
&&\qquad= \bigl(1 + O\bigl(e^{-\gamma c k} + g(y,z) + c k \Upsilon_n\bigr)\bigr) \pi(z)
\end{eqnarray*}
as desired.
\end{pf*}

In the following lemma, we will show that the difference in the
probability that a random walk hits points $y,z$ when started from $x$
before time $\Gamma_n = c_0 ({\log}|V_n|) T_\mix$ is essentially
determined by the corresponding difference except up to time $c T_\mix
$. The reason for the cancellation is that the previous lemma implies
that conditional on not hitting a given point up to time $c T_\mix$,
the walk is well mixed and has long forgotten its starting point.
Recall that $f_c(x,y) = \p_x[ \tau(y) \leq cT_\mix]$ (we have
suppressed $n$).
%
%
\begin{lemma}
\label{lemmixingcancel}
There exists $\gamma> 0$ such that for all $c \geq2$,
\begin{eqnarray*}
&&\p_x[ \tau(y) \leq\Gamma_n] - \p_x[ \tau(z) \leq\Gamma_n] \\
&&\qquad= f_c(x,y)
- f_c(x,z)
 + O\bigl( e^{-\gamma c} \Upsilon_n + \Delta
_n[g(x,y) + g(x,z)]\bigr).
\end{eqnarray*}
\end{lemma}
\begin{pf}
We observe
\begin{eqnarray*}
&&\p_x[ \tau(y) \leq\Gamma_n] \\
&&\qquad= f_c(x,y)
+\sum_{k} \p_x[ c T_k < \tau(y) \leq c T_{k+1} | \tau(y) > c T_k]
\bigl(1-\p
_x[ \tau(y) \leq c T_k]\bigr),
\end{eqnarray*}
where, here and throughout the rest of this proof, the summation over
$k$ is from $1$ to ${\frac{c_0}{c} \log}|V_n|$. We note that
\begin{eqnarray*}
&& \p_x[c T_k < \tau(y) \leq c T_{k+1} | \tau(y) > c T_k]\\
&&\qquad= \sum_w \p_x[ c T_k < \tau(y) \leq c T_{k+1}| \tau(y) > c T_k,
X(cT_k) = w] \\
&&\qquad\quad\hspace*{13pt}{}\times\p_x[X(cT_k) = w|\tau(y) > c T_k]\\
&&\qquad= \sum_w \p_w[ \tau(y) \leq c T_{1}] \p_x[X(cT_k) = w|\tau(y) > c T_k].
\end{eqnarray*}
As the previous lemma is applicable for such choices of $k$ and using
$\p_w[\tau(y) \leq c T_1] \leq O(g(y,w) + c \Upsilon_n)$, we can
rewrite the expression above as
\[
\p_\pi[\tau(y) \leq c T_1] + O\biggl(\sum_{w \neq y} g(y,w)
\bigl(e^{-\gamma
c k} + c k \Upsilon_n + g(y,w)\bigr) \pi(w) \biggr).
\]
Performing the summation over $w$, we see that the latter term is of order
%
%
\begin{equation}
\label{eqne1}
O(\Upsilon_n e^{-\gamma c k} + c k \Upsilon_n^2 + \CS_n |V_n|^{-1}).
\end{equation}
Recall from part~\ref{assumpmaingreensquarebound} of Assumption
\ref{assumpmain} that $\CS_n = o(T_\mix/ {\log}|V_n|)$, hence
$({\log}
|V_n|) \CS_n |V_n|^{-1} = o(\Upsilon_n)$. Consequently, summing
(\ref{eqne1}) over $k$ from $1$ to $\frac{c_0}{c}\times{\log}|V_n|$ gives an
error of
\[
O( \Upsilon_n e^{-\gamma c} + \Delta_n^2).
\]
By part~\ref{assumpmainmixingbound} of Assumption \ref
{assumpmain} it is clear that $\Delta_n^2 = o(\Upsilon_n)$, hence the
former is of order $O(\Upsilon_n e^{-\gamma c})$.
This leaves
\begin{eqnarray*}
&&\sum_k \p_x[ c T_k < \tau(y) \leq c T_{k+1} | \tau(y) > c T_k] \p_x[
\tau(y) \leq c T_k]\\
&&\qquad= O\biggl( \sum_k \sum_{z} \p_z[ \tau(y) \leq c T_1] \pi(z)
\bigl(g(x,y) + c
k \Upsilon_n\bigr) \biggr).
\end{eqnarray*}
Here, we used the previous lemma to get the crude estimate $\p_x[X(c
T_k) = z| \tau(y) > c T_k] \leq C\pi(z)$ for some $C > 0$. Summing
everything up gives us an error of order $O( \Delta_n g(x,y) + \Delta
_n^2)$. We also have another contribution of $O(\Delta_n g(x,z) +
\Delta
_n^2)$ coming from the corresponding estimate of $\p_x[\tau(z) \leq
\Gamma_n]$. Therefore our total error is $O(\Upsilon_n e^{-\gamma c} +
\Delta_n[g(x,y) + g(x,z)])$, which proves the lemma.
\end{pf}

\section{The Radon--Nikodym derivative}
\label{secrn}

Recall
\[
\wt{\pi}(z;x,y) = \p[ X_2(\tau_1(x,y)) = z | H(x,y)].
\]
The purpose of this section is to prove the following estimate of the
Radon--Nikodym derivative of $\wt{\pi}(z;x,y)$ with respect to $\pi
(z)$. Recall\vadjust{\goodbreak} again $f_c(x,y) = \p_x[ \tau(y) \leq c T_\mix]$ and
$\ol
{f}_c = \sum_y f_c(x,y) \pi(y)$ (we are omitting the dependence on~$n$).
%
%
\begin{theorem}
\label{thmrnestimate}
There exists a constant $\gamma> 0$ so that for all $c \geq2$ and $x
\neq y$, we have
\begin{eqnarray*}
\frac{\wt{\pi}(z;x,y)}{\pi(z)}
&=& 1 + \bigl(1 + O(\Delta_n) \bigr) \bigl( 2\ol{f}_c - f_c(x,z) -
f_c(y,z) \bigr)\\
&&{} + O(e^{-\gamma c} \Upsilon_n ) +
O\bigl([ g(x,z) + g(y,z) + \Delta_n][g(x,y) + \Delta_n]\bigr).
\end{eqnarray*}
In particular,
%
%
\begin{equation}
\label{eqnradonsimplebound}
\frac{\wt{\pi}(z;x,y)}{\pi(z)} = 1 + O\bigl(g(x,y) + g(y,z) + g(x,z)\bigr).
\end{equation}
\end{theorem}

\begin{figure}

\includegraphics{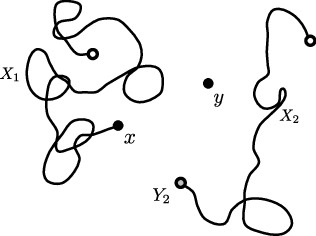}

\caption{In Theorem \protect\ref{thmrnestimate}, we give
a precise estimate of the Radon--Nikodym derivative of the law of $Y_2
= X_2(\tau_1(x,y))$ with respect to the uniform measure on $V_n$
conditional on the event that $H(x,y) = \{\tau_1(x) < \tau_1(y)
\wedge
\tau_2(x,y)\}$, that is, that the first point in $\{x,y\}$ hit by
$X_1,X_2$ is $x$ by~$X_1$. The open circles indicate the starting
points of $X_1,X_2$ and the shaded circle is $Y_2$.}\label{fig1}
\end{figure}

\begin{figure}

\includegraphics{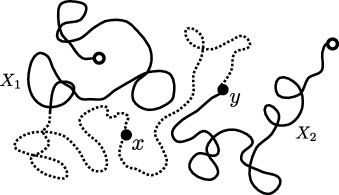}

\caption{Illustration of the event $A(x,y)$. The solid lines are used
to indicates the parts of $X_1,X_2$ up to the time $\tau_2(x,y)$, while
the dashed line is used for the part of $X_1$ after $\tau_2(x,y)$. We
have not indicated the part of $X_2$ after $\tau_2(x,y)$. Note that we
have $X_2$ hitting $y$ first, but $A(x,y)$ allows for $X_2$ to hit $x$
first as well. On the other hand, $A(x,y)$ requires that $X_1$ does in
fact hit $x$ before $y$.}\label{fig2}
\end{figure}

The setup for Theorem~\ref{thmrnestimate} is illustrated in Figure~\ref{fig1}.
%
%
Let $Y_2 = X_2(\tau_1(x,y))$. The idea of the proof is to observe that
\[
\wt{\pi}(z;x,y)
= \p[ Y_2 = z | H(x,y)]
= \frac{\p[ H(x,y) | Y_2 = z] \pi(z)}{\p[ H(x,y)]},
\]
where we used $\p[ Y_2 = z] = \pi(z)$ as $X_1,X_2$ are independent and
the initial distribution of $X_2$ is stationary, then estimate the
effect of conditioning on $\{Y_2 = z\}$ on the probability of $H(x,y)$.
We will divide the proof into three lemmas. The first step in the proof
is to express $\wt{\pi}(\cdot;x,y)/\pi$ in terms of the event
\begin{eqnarray*}
A(x,y)
&=& \{ \tau_2(x,y) > \tau_1(x,y) - \Gamma_n, G_1(x,y)\} \setminus
H(x,y)\\
&=& \{\tau_1(x,y) \geq\tau_2(x,y) > \tau_1(x,y) - \Gamma_n,
G_1(x,y) \}.
\end{eqnarray*}
The event $A(x, y)$ is illustrated in Figure~\ref{fig2}.
Note that it is a slight abuse of notation to insert $G_1(x,y)$ into
the braces defining $A(x,y)$ since $G_1(x,y)$ is itself an event. We
%
%
will do this a number of times in the following lemma in order to
lighten the notation.

%
\begin{lemma}
\label{lemradonbound}
Uniformly in $x,y,z,n$,
%
%
\begin{equation}\quad
\frac{\wt{\pi}(z;x,y)}{\pi(z)}
=1 + \frac{\p[A(x,y)] - \p[A(x,y) | Y_2 = z]}{\p[H(x,y)]} +
O(|V_n|^{-100}).
\end{equation}
\end{lemma}
\begin{pf}
Letting $R(x,y) = \{\tau_2(x,y) > \tau_1(x,y) - \Gamma_n, G_1(x,y)\}
$, observe
\[
\p[H(x,y) | Y_2 = z]
= \p[R(x,y)| Y_2 = z]
- \p[A(x,y)| Y_2 = z].
\]
We will now manipulate the first term on the right-hand side. Let $\wt
{R}(x,y) = G_1(x,y) \setminus R(x,y)$. We have
%
%
\begin{equation}\label{eqnrnbreakdown1}
\p[R(x,y) | Y_2 = z]
= \p[G_1(x,y)| Y_2 = z] - \p[\wt{R}(x,y) | Y_2 = z]
\end{equation}
and, since $Y_2 \sim\pi$, Bayes' rule implies
\[
\p[\wt{R}(x,y) | Y_2 = z]
= \frac{1}{\pi(z)}\p[Y_2 = z| \wt{R}(x,y)]
\p[\wt{R}(x,y)].
\]
Since the conditional probability on the right-hand side involves
conditioning on the behavior of $X_2$ before $\tau_1(x,y) - \Gamma_n$;
mixing considerations imply that this is equal to
%
%
\begin{equation}\label{eqnrnbreakdown2}
[1+O(|V_n|^{-\ol{\gamma} c_0})]\p[ \wt{R}(x,y)]
= \p[\wt{R}(x,y)] + O(|V_n|^{-\ol{\gamma} c_0}) .
\end{equation}
As $\p[ G_1(x,y) | Y_2 = z] = \p[ G_1(x,y)]$, combining (\ref
{eqnrnbreakdown1}) with (\ref{eqnrnbreakdown2}) we thus have
\begin{eqnarray*}
\p[ H(x,y) | Y_2 = z]
&=& \p[R(x,y)]
- \p[ A(x,y) | Y_2 = z] + O(|V_n|^{-\ol{\gamma} c_0})\\
&=& \p[H(x,y)] + \p[A(x,y)] - \p[ A(x,y) | Y_2 = z]\\
&&{} + O(|V_n|^{-\ol
{\gamma} c_0}).
\end{eqnarray*}
Assume that $\ol{\gamma} c_0 > 100$.
Putting everything together, we see that
\[
\frac{\wt{\pi}(z;x,y)}{\pi(z)} = \frac{\p[H(x,y)] + \p[A(x,y)] -
\p[A(x,y)|Y_2 = z]}{\p[H(x,y)]} + O(|V_n|^{-100}),
\]
uniformly in $x,y,z,n$.
\end{pf}

Note that if $V_n = \Z_n^d$ for $d \geq3$, then $\p[ G_{1}(x,y),
G_{12}(x,y)] = \p[H(x,y)]$ does not change when $x$ is swapped with $y$
nor when $1$ is swapped with $2$ and hence is equal to $\frac{1}{4}$
up to negligible error (it is not exactly $\frac{1}{4}$ since it could
be that $X_1$ hits $x$ at the same time $X_2$ hits either $x$ or $y$,
though this is a rare event). This holds more generally if for every
$x,y \in V_n$ distinct there exists an automorphism $\varphi$ of $G_n$
such that $\varphi(x) = y$ and $\varphi(y) = x$. The weaker hypothesis
of vertex transitivity implies that we can always find an automorphism
$\varphi$ of $G_n$ such that $\varphi(x) = y$ but not necessarily so
that $\varphi(y) = x$ as well. Nevertheless, it is still true in this
case that $\p[H(x,y)] \approx\frac{1}{4}$.

%
\begin{lemma}
\label{lemsymmetry}
If $x \neq y$, we have that
\[
\p[H(x,y)] = \frac{1}{4} + o\biggl( \frac{\Upsilon_n}{ {\log}|V_n|}
\biggr) + O\biggl( \frac{1}{|V_n|} \biggr).
\]
\end{lemma}
\begin{pf}
Let $\wt{A}(x,y) = \{ \tau_1(x,y) \geq\tau_2(x,y) > \tau_1(x,y) -
\Gamma_n\}$ and $\mu(z;x,y) = \p[ Y_2 = z| \tau_1(x,y) \leq\tau
_2(x,y)]$. Using exactly the same proof as the previous lemma, we have
\[
\frac{\mu(z;x,y)}{\pi(z)} = 1 + O\bigl(\p[\wt{A}(x,y)] + \p[\wt{A}(x,y)|Y_2
= z] + |V_n|^{-100}\bigr).
\]
Using a time-reversal in the first step and a union bound in the
second, we have that
\[
\p[\wt{A}(x,y) | Y_2 = z] \leq\p_z[ \tau_2(x,y) \leq\Gamma_n] =
O\bigl(\Delta_n + g(x,z) + g(y,z)\bigr)
\]
(and similarly for $\p[\wt{A}(x,y)]$). Consequently,
\begin{eqnarray*}
\sum_z g(x,z) \mu(z;x,y)
&=& \Upsilon_n + \frac{1}{|V_n|}\sum_z g(x,z) O\bigl( \Delta_n + g(x,z)
+ g(y,z) \bigr)\\
&=& \Upsilon_n + \frac{1}{|V_n|}\sum_z O\bigl(g(x,z) \Delta_n +
g^2(x,z) + g^2(y,z) \bigr)\\
&=& \Upsilon_n + \frac{1}{|V_n|}O(1+\CS_n + \Delta_n T_\mix).
\end{eqnarray*}
By parts~\ref{assumpmainmixingbound} and \ref
{assumpmaingreensquarebound} of Assumption~\ref{assumpmain}, we
have that $\CS_n + \Delta_n T_\mix= o(T_\mix/({\log}|V_n|))$.
Consequently, the above is equal to
%
%
\begin{equation}
\label{eqnapproxunifgreen}
\Upsilon_n +
o\biggl(\frac{\Upsilon_n}{{\log}|V_n|} \biggr).\vadjust{\goodbreak}
\end{equation}
Let $p_x = \p[G_1(x,y) , G_{12}(x,y)]$ and $p_y = \p[ G_1(y,x),
G_{12}(x,y)]$. Note that
%
%
\begin{equation}
\label{eqnpxpy} p_x+p_y = \p[ G_{12}(x,y)] = \frac{1}{2} + O\biggl(
\frac{1}{|V_n|} \biggr)
\end{equation}
since $\p[\tau_1(z) = \tau_2(w)] \leq\p[X_2(\tau_1(z)) = w] =
|V_n|^{-1}$ for any $z,w \in V_n$. Define stopping times as follows. Let
\[
\tau_1 = \min\bigl\{t \geq0\dvtx X_1(t) \in\{x,y\} \mbox{ or } X_2(t) \in
\{
x,y\}\bigr\} = \tau(x,y).
\]
For $j \geq1$, inductively set
\[
\tau_{j+1} = \min\bigl\{ t \geq\tau_j + T_\mix+ 1\dvtx X_1(t) \in\{x,y\}
\mbox{ or } X_2(t) \in\{x,y\}\bigr\}.
\]
Let $\CT_{j,z} = \sum_{t = \tau_j}^{\tau_j+T_\mix} \one_{\{X_1(t) =
z\}
}$, and, for $E \subseteq V_n$, set $A_{ij}(E) = \{ X_i(\tau_j) \in E\}
$. Note that
the average amount of time spent at $x$ by $X_1$ through time $\tau_k +
T_\mix$ is given by the expression
\[
\frac{1}{\tau_k+T_\mix} \sum_{j=1}^k \bigl( \one_{A_{1j}(x)} \one
_{A_{2j}^c(x,y)} \CT_{j,x} +
\one_{A_{1j}(y)} \one_{A_{2j}^c(x,y)} \CT_{j,x} +
\one_{A_{2j}(x,y)} \CT_{j,x} \bigr).
\]
It is not difficult to see that the above quantity converges to $\pi
(x)$ as $k \to\infty$. We can also define a similar quantity but
replacing $\CT_{j,x}$ with $\CT_{j,y}$; this will converge to $\pi(y)$
as $k \to\infty$. Taking the ratio of these two quantities, we arrive at
\[
1 = \lim_{k \to\infty}
\frac{ ({1}/{k}) \sum_{j=1}^k ( \one_{A_{1j}(x)} \one
_{A_{2j}^c(x,y)} \CT_{j,x} +
\one_{A_{1j}(y)} \one_{A_{2j}^c(x,y)} \CT_{j,x} +
\one_{A_{2j}(x,y)} \CT_{j,x} )}
{({1}/{k}) \sum_{j=1}^k ( \one_{A_{1j}(x)} \one_{A_{2j}^c(x,y)}
\CT_{j,y} +
\one_{A_{1j}(y)} \one_{A_{2j}^c(x,y)} \CT_{j,y} +
\one_{A_{2j}(x,y)} \CT_{j,y} )}
\]
since $\pi(x) = \pi(y)$.
It is not difficult to see that, almost surely,
\begin{eqnarray*}
\lim_{k \to\infty} \frac{1}{k} \sum_{j=1}^k \one_{A_{1j}(x)}
\one
_{A_{2j}^c(x,y)} \CT_{j,x} &=& p_x g(x,x),\\
\lim_{k \to\infty} \frac{1}{k} \sum_{j=1}^k \one_{A_{1j}(y)}
\one
_{A_{2j}^c(x,y)} \CT_{j,x} &=& p_y g(y,x),\\
\lim_{k \to\infty} \frac{1}{k} \sum_{j=1}^k \one_{A_{2j}(x,y)}
\CT
_{j,x} &=& q_{xy} \sum_{z} g(z,x) \mu(z;x,y),
\end{eqnarray*}
where $q_{xy} = 1-p_x-p_y$. Analogous formulae hold for the terms in
the denominator. Combining this with (\ref{eqnapproxunifgreen}), we
thus have
\begin{eqnarray*}
1 & = &\frac{p_x g(x,x) + p_y g(y,x) + q_{xy} \sum_{z} g(z,x) \mu
(z;x,y)}{p_y g(y,y) + p_x g(x,y) + q_{xy} \sum_z g(z,y) \mu(z;x,y)}\\
&=&\frac{p_x g(x,x) + p_y g(y,x) + q_{xy} \Upsilon_n}{p_y g(y,y) + p_x
g(x,y) + q_{xy} \Upsilon_n} + o\biggl( \frac{\Upsilon_n}{\log
|V_n|}\biggr).
\end{eqnarray*}
Rearranging and using that $g(x,y) = g(y,x)$ and $g(x,x) = g(y,y)$,
this implies that
\[
p_x = p_y + o\biggl( \frac{\Upsilon_n}{ {\log}|V_n|} \biggr).
\]
Combining this with (\ref{eqnpxpy}) proves the lemma.
\end{pf}

In order to complete the proof of Theorem~\ref{thmrnestimate} we
need to estimate $\p[A(x,y)| Y_2 = z]$, which is the purpose of the
following lemma. Though the proof will be computationally intensive,
the basic idea is fairly simple. The main goal is to eliminate the
conditioning on $Y_2 = z$. The first step is to perform a time
reversal, which converts the terminal condition to an initial condition
at the cost of making the event whose probability we are to compute a
bit more complicated. The latter is easily mitigated, however, since
the event can be greatly simplified at the cost of negligible error.
%
%
\begin{lemma}
\label{lemabound}
There exists $\gamma> 0$ so that for all $c \geq2$ we have
\begin{eqnarray*}
\p[ A(x,y) | Y_2 = z]
&=& \bigl( \tfrac{1}{4} + O(\Delta_n) \bigr) [ f_c(x,z) + f_c(y,z)
] + E_c(x,y) \\
&&{}+ O\bigl( e^{-\gamma c} \Upsilon_n + [g(x,z) + g(y,z) + \Delta
_n][g(x,y) + \Delta_n] \bigr),
\end{eqnarray*}
where $E_c(x,y)$ is some constant which does not depend on $z$.
\end{lemma}

Note that the lemma implies
\[
\p[ A(x,y)] - \p[ A(x,y) | Y_2 = z]
= O\bigl(g(x,y) + g(y,z) + g(x,z) + e^{-\gamma c} \Upsilon_n\bigr).
\]
\begin{pf*}{Proof of Lemma~\ref{lemabound}}
Let
\[
B(x,y) = \bigl\{ X_2(t) \notin\{x,y\} \mbox{ for all } t \in(\Gamma_n,
\tau
_1(x,y)], G_1(x,y)\bigr\},
\]
and let $\p_{\pi,z}$ be the law of $(X_1,X_2)$ where $X_1(0) \sim\pi$
and $X_2(0) = z$. We compute
\begin{eqnarray*}
&&
\p[ A(x,y) | Y_2 = z]\\
&&\qquad= \frac{1}{\pi(z)} \p[ A(x,y), Y_2 = z]\\
&&\qquad= \sum_w \p_{\pi,w}[ \tau_1(x,y) \geq\tau_2(x,y) > \tau_1(x,y) -
\Gamma_n, G_{1}(x,y), Y_2 = z].
\end{eqnarray*}
By reversing the time of $X_2$ (but not $X_1$), we see that this is
equal to
%
%
\begin{eqnarray}\label{eqnreversebound}
&& \sum_w \p_{\pi,z}[ \tau_2(x,y) \leq\Gamma_n \wedge\tau_1(x,y),
B(x,y), Y_2 = w] \nonumber\\[-8pt]\\[-8pt]
&&\qquad= \p_{\pi,z}[ \tau_2(x,y) \leq\Gamma_n \wedge\tau_1(x,y), B(x,y)].
\nonumber
\end{eqnarray}
We will now work toward approximating this event with a simpler event.
We begin by eliminating the ``minimum'' operation using the observation
that it is unlikely for both $X_1, X_2$ to hit $\{x,y\}$ quickly.
Indeed, as
\[
\p_{\pi,z}[\tau_1(x,y) \leq\Gamma_n, \tau_2(x,y) \leq\Gamma_n]
= O\bigl( [ g(x,z) + g(y,z) + \Delta_n] \Delta_n \bigr),
\]
we see by setting $\wt{B}(x,y) = B(x,y) \cap\{\tau_1(x,y) > \Gamma
_n\}
$ that (\ref{eqnreversebound}) is equal to
\[
\p_{\pi,z}[ \tau_2(x,y) \leq\Gamma_n, \wt{B}(x,y)]
+ O\bigl( [ g(x,z) + g(y,z) + \Delta_n] \Delta_n\bigr).
\]
We would now like to eliminate the dependence of the probability on
$z$, the starting point of $X_2$. We accomplish this by considering two
possible cases. Either $X_2$ hits $x$ or $y$ within some multiple of
the mixing time or it does not. Conditional on the latter, the walk
will have mixed, so the relevant probability does not depend on $z$. We
implement this strategy as follows:
\begin{eqnarray*}
&&\p_{\pi,z}[ \tau_2(x,y) \leq\Gamma_n, \wt{B}(x,y)]\\
&&\qquad= \p_{\pi,z}[ \tau_2(x,y) < c T_\mix, \wt{B}(x,y)] \\
&&\qquad\quad{}+ \p_{\pi,z}[ \tau_2(x,y) \leq\Gamma_n, \wt{B}(x,y) |
\tau_2(x,y) \geq c T_\mix]\\
&&\qquad\quad\hspace*{11pt}{}\times \bigl(1-\p_z[ \tau_2(x,y) < c T_\mix]\bigr).
\end{eqnarray*}
Using the same proof as Lemma~\ref{lemconditionnothitrn}, except
in the case that the random walk is conditioned not to hit two points
rather than just one, implies $\mu(w;x,y,z) = \p_z[ X_2(c T_\mix) =
w |
\tau_2(x,y) \geq c T_\mix] \leq C\pi(w)$ for some constant $C > 0$.
Consequently,
\begin{eqnarray*}
&& \p_{\pi,z}[ \tau_2(x,y) \leq\Gamma_n, \wt{B}(x,y) | \tau
_2(x,y) \geq
c T_\mix] \p_z[ \tau_2(x,y) < c T_\mix] \\
&&\qquad\leq C\p_\pi[ \tau_2(x,y) \leq\Gamma_n] \p_z[ \tau_2(x,y) < c
T_\mix]\\
&&\qquad= O\bigl( [ g(x,z) + g(y,z) + c \Upsilon_n] \Delta_n\bigr).
\end{eqnarray*}

We are left with two terms to estimate
%
%
\begin{eqnarray}
\label{eqnhitearly}
&\p_{\pi,z}[ \tau_2(x,y) < c T_\mix, \wt{B}(x,y)],&
\\
\label{eqnhitlate}
&\p_{\pi,z}[ \tau_2(x,y) \leq\Gamma_n, \wt{B}(x,y) | \tau_2(x,y)
\geq
c T_\mix] .&
\end{eqnarray}
We will first deal with (\ref{eqnhitearly}) which, using the
independence of $\tau_1(x,y)$ and $\tau_2(x,y)$, we can rewrite as
\begin{eqnarray*}
&&\p_{\pi,z}[B(x,y)|\tau_1(x,y) > \Gamma_n , \tau_2(x,y) < c T_\mix
] \\
&&\qquad{}\times\p
_{z}[\tau_2(x,y) < c T_\mix]
\p[\tau_1(x,y) > \Gamma_n].
\end{eqnarray*}
Since $B(x,y)$ depends on $X_2(t)$ only for $t \geq\Gamma_n$, from
mixing considerations it is easy to see that
\begin{eqnarray*}
&& \p_{\pi,z}[B(x,y)|\tau_1(x,y) > \Gamma_n , \tau_2(x,y) < c
T_\mix]\\
&&\qquad= \p[B(x,y)|\tau_1(x,y) > \Gamma_n] + O(|V_n|^{-100}).
\end{eqnarray*}
Consequently, (\ref{eqnhitearly}) is equal to
\[
\p_z[ \tau_2(x,y) < c T_\mix] \p[\wt{B}(x,y)] + O(|V_n|^{-100}).
\]
Note that
\[
\wt{B}(x,y) = \bigl(H(x,y) \cap\{ \tau(x,y) > \Gamma_n\}\bigr)
\cup\bigl(\wt{B}(x,y)
\cap\{\tau(x,y) \leq\Gamma_n\}\bigr).
\]
Using $\p[ \tau(x,y) \leq\Gamma_n] = O(\Delta_n)$, the previous lemma
thus implies
\[
\p[\wt{B}(x,y)] = \p[ H(x,y) ] + O(\Delta_n) = \tfrac{1}{4} +
O(\Delta_n).
\]
Observe
\begin{eqnarray*}
&& \p_z[ \tau_2(x) < c T_\mix, \tau_2(y) < c T_\mix]\\
&&\qquad= \p_z[ \tau_2(x) < \tau_2(y) < c T_\mix] + \p_z[\tau_2(y) <
\tau_2(x)
< c T_\mix]\\
&&\qquad= O\bigl( [g(x,z) + g(y,z) + c \Upsilon_n][g(x,y) + c \Upsilon_n]\bigr).
\end{eqnarray*}
Consequently,
\begin{eqnarray*}
&&\p_z[ \tau_2(x,y) < c T_\mix]\\
&&\qquad= f_c(x,z) + f_c(y,z) - \p_z[ \tau_2(x) < c T_\mix, \tau_2(y) < c
T_\mix]\\
&&\qquad= f_c(x,z) +f_c(y,z) + O\bigl([g(x,z) + g(y,z) + c \Upsilon_n][g(x,y) + c
\Upsilon_n]\bigr).
\end{eqnarray*}
%
Arguing as in the proof of
Lemma~\ref{lemmixingcancel}, we can estimate (\ref{eqnhitlate})
as follows:
\begin{eqnarray*}
&&\p_{\pi,z}[ \tau_2(x,y) \leq\Gamma_n, \wt{B}(x,y) | \tau_2(x,y)
\geq
c T_\mix]\\
&&\qquad= \p[ \tau_2(x,y) \leq\Gamma_n, \wt{B}(x,y) | \tau_2(x,y) \geq c
T_\mix]\\
&&\qquad\quad{} + O\bigl(e^{-\gamma c} \Upsilon_n + \bigl(g(x,z) + g(y,z)\bigr) \Delta_n\bigr).
\end{eqnarray*}
Taking $E_c(x,y) = \p[ \tau_2(x,y) \leq\Gamma_n, \wt{B}(x,y) |
\tau
_2(x,y) \geq c T_\mix]$ and noting that $c \Upsilon_n = O(\Delta_n)$
finishes the proof of the lemma.
\end{pf*}

By combining the three lemmas, we can now complete the proof of
Theorem~\ref{thmrnestimate}.
\begin{pf*}{Proof of Theorem~\ref{thmrnestimate}}
Recalling that $\p[Y_2 = z] = \pi(z)$, observe
\[
\p[ A(x,y)]
= \sum_z \p[ A(x,y) | Y_2 = z] \p[Y_2 = z]
= \sum_z \p[ A(x,y) | Y_2 = z] \pi(z).
\]
Hence by Lemma~\ref{lemabound}, we have
\[
\p[ A(x,y)]
= \bigl( \tfrac{1}{4} + O(\Delta_n) \bigr)(2\ol{f}_c) + E_c(x,y) +
O\bigl( e^{-\gamma c} \Upsilon_n + \Delta_n[g(x,y) + \Delta_n] \bigr).
\]
Here, we used that
\[
\sum_z \bigl(g(x,z) + g(y,z)\bigr) \pi(z) = O(\Upsilon_n) = O(\Delta_n).
\]
In particular,
\begin{eqnarray*}
&&\p[A(x,y)] - \p[A(x,y) | Y_2 = z]\\
&&\qquad= \bigl( \tfrac{1}{4} + O(\Delta_n) \bigr) [ 2\ol{f}_c -
f_c(x,z) -
f_c(y,z) ] \\
&&\qquad\quad{} + O\bigl( e^{-\gamma c} \Upsilon_n + [g(x,z) + g(y,z) + \Delta
_n][g(x,y) + \Delta_n] \bigr).
\end{eqnarray*}
Inserting this expression along with the estimate of $\p[H(x,y)]$ from
Lem\-ma~\ref{lemsymmetry} into into the equation for $\wt{\pi}(z;x,y)/
\pi(z)$ from Lemma~\ref{lemradonbound} gives the theorem.
\end{pf*}

\section{The variance}
\label{secvar}

We will complete the proof of Theorem~\ref{thmmain} in this section.
The general theme is to eliminate asymmetry wherever possible. We first
apply this idea by considering
\[
\CB= \sum_x \one_{G_{12}(x)} - \sum_x \one_{G_{21}(x)}
\]
in place of $|\CA_1|$. In addition to being symmetric in $X_1, X_2$,
note that $\CB$ also differs from $|\CA_1|$ in that we have eliminated
those sites whose mark is determined by the flip of a fair coin. These,
however, do not make a significant contribution to the variance since
it is a rare event that both walks hit a particular point for the first
time simultaneously.
In particular, we will show in Lemma~\ref{lemsymm} that $\var(\CB)
\approx4 \var(|\CA_1|)$, up to negligible error. Consequently, to
prove Theorem~\ref{thmmain} it suffices to show
\[
\var(\CB) = \sum_{x,y} \bigl(f_c(x,y) - \ol{f}_c\bigr)^2
+ O(e^{-\gamma c}
(T_\mix)^2).
\]
It is convenient to work with $\CB$ as the expansion of its variance
takes on the following form:
%
%
\begin{eqnarray}
\var(\CB)
&=& 2\sum_{x,y} \bigl( \p[ G_{12}(x), G_{12}(y)] - \p[ G_{12}(x),
G_{21}(y)] \bigr) \nonumber\\
\label{eqnvaroffdiag}
&=& 4\sum_{x \neq y \neq z}\bigl( \p_{x,z}[ G_{12}(y)] - \p_{x,z}[
G_{21}(y)] \bigr) \wt{\pi}(z;x,y) \p[H(x,y)]
\\
\label{eqnvardiag}
&&{}+ 2 \sum_{x} \p[G_{12}(x)].
\end{eqnarray}
The reason the summation in (\ref{eqnvaroffdiag}) is over $x \neq y
\neq z$ is $\p[H(x,x)] = 0$ and $\wt{\pi}(z;x,y) = 0$ if $z \in\{
x,y\}
$; the summation in (\ref{eqnvardiag}) contains the diagonal terms.
We will now focus on (\ref{eqnvaroffdiag}) and handle (\ref
{eqnvardiag}) at the end of the section. Applying Lemma \ref
{lemsymmetry}, we can rewrite (\ref{eqnvaroffdiag}) as
%
%
\begin{eqnarray}
\label{eqnvarianceformula}
&& \sum_{x \neq y \neq z} \bigl( \p_{x,z}[ G_{12}(y)] - \p_{x,z}[
G_{21}(y)] \bigr) \wt{\pi}(z;x,y)  \\
\label{eqnvarianceneg}
&&\qquad{}+ \sum_{x \neq y \neq z}|\p_{x,z}[ G_{12}(y)] - \p_{x,z}[G_{21}(y)]|
\nonumber\\[-8pt]\\[-8pt]
&&\hspace*{11pt}\qquad{}\times\biggl[ o\biggl( \frac{T_\mix}{{|V_n|^2 \log}|V_n|}\biggr) + O\biggl(
\frac
{1}{|V_n|^2}\biggr) \biggr] .\nonumber
\end{eqnarray}
We will show at the end of this section that (\ref{eqnvarianceneg})
is negligible.
Note that
%
%
\begin{eqnarray}\label{eqntimedecomp}
\p_{x,z}[ G_{ij}(y)]
&=& \p_{x,z}[ \tau_i(y) < \tau_j(y) \leq\Gamma_n] \nonumber\\
&&{}+ \p_{x,z}[ \tau
_i(y) \leq\Gamma_n, \tau_j(y) > \Gamma_n]
+ \p_{x,z}[ \Gamma_n < \tau_i(y) < \tau_j(y)] \hspace*{-25pt}\\
&\equiv& A + B + C.
\nonumber
\end{eqnarray}
In Section~\ref{subsecconditioning}, we break the sum in (\ref
{eqnvarianceformula}) into three different cases based on the time
decomposition in (\ref{eqntimedecomp}) and bound each in a given
lemma. It will turn out that the contributions to the variance coming
from the terms corresponding to $A$ and $C$ are negligible (Lemmas \ref
{lema} and~\ref{lemc}). The reason for the former is that it
is unlikely for both $X_1$ and $X_2$ to hit $y$ quickly and the latter
follows as, conditional on having not hit $y$ by time $\Gamma_n$, both
walks have long forgotten their initial conditions and are well mixed.
This leaves $B$, which, along with the diagonal, dominates the
variance. Its asymptotics will be computed (Lemma~\ref{lemb}) by
reducing the estimate to a computation involving $\wt{\pi}(z;x,y)$,
whose Radon--Nikodym derivative with respect to the uniform measure has
already been estimated precisely in Theorem~\ref{thmrnestimate}.

\subsection{Symmetrization}
%
%
\begin{lemma}
\label{lemsymm}
We have
\[
\var(\CB) = 4 \var(|\CA_1|) + O\bigl(\sqrt{T_\mix\var(|\CA
_1|)} +
T_\mix\bigr) \qquad\mbox{as } n \to\infty.
\]
\end{lemma}
\begin{pf}
Let $(\xi_{n}(x)\dvtx x \in V_n)$ be i.i.d. random variables independent
of $X_1,X_2$ with $\p[\xi_n(x) = 1] = \p[\xi_n(x) = 2] = \frac{1}{2}$,
and let $A(x,i) = \{ \tau_1(x) = \tau_2(x), \xi_n(x) = i\}$. By
definition,
\begin{eqnarray*}
\var(\CB)
&=& \var\biggl( |\CA_1| - ( |V_n| - |\CA_1| ) - \sum_{x}
\bigl(
\one_{A(x,1)} - \one_{A(x,2)} \bigr) \biggr)\\
&=& \var\biggl( 2|\CA_1| + \sum_{x} \bigl(\one_{A(x,2)} - \one_{A(x,1)}
\bigr) \biggr).
\end{eqnarray*}
Observe
\[
\E\biggl( \sum_{x} \one_{A(x,1)} \biggr)^2
\leq\sum_{x,y} \p[ \tau_1(x) = \tau_2(x), \tau_1(y) = \tau_2(y)].
\]
By the strong Markov property and independence of $X_1,X_2$, the above
is bounded by twice
\begin{eqnarray*}
&&\sum_{x,y} \p_{x,x}[ \tau_1(y) = \tau_2(y)] \p[ \tau_1(x) = \tau
_2(x)] \\
&&\qquad\leq\sum_{x,y} \sum_{t} \bigl(\p_{x}[ \tau(y) = t] \bigr)^2 \p[
X_2(\tau
_1(x)) = x].
\end{eqnarray*}
Using that $\p_x[\tau(y) = t] \leq\p_x[X(t) = y]$ and $X_2(\tau_1(x))
\sim\pi$ when $X(0) \sim\pi$, we have the further bound
%
%
\begin{equation}\label{eqnhitsametime}
\frac{1}{|V_n|}\sum_{x,y} \Biggl( \sum_{t=0}^{4T_\mix} \p_{x}[
X(t) = y]
+ \sum_{t > 4T_\mix} \bigl(\p_x[\tau(y) = t]\bigr)^2 \Biggr).
\end{equation}
Summing the first term over $x,y,t$ plainly yields $4T_\mix$. For the
second term, note there exists $C > 0$ so that for $t > 4T_\mix$, we have
\[
\p_x[\tau(y) = t] \leq\p_x[ X(t) = y] \leq\frac{C}{|V_n|}
\]
hence
\[
\sum_{t > 4 T_\mix} \bigl(\p_x[\tau(y) = t]\bigr)^2 \leq\frac{C}{|V_n|}
\sum_{t
> 4 T_\mix} \p_x[\tau(y) = t] \leq\frac{C}{|V_n|}.
\]
Therefore the second term in the summation in (\ref{eqnhitsametime})
is $O(1)$. The lemma now follows from Cauchy--Schwarz.
\end{pf}

\subsection{Time decomposition}
\label{subsecconditioning}

We begin by estimating the part of (\ref{eqnvarianceformula})
corresponding to ``$B$'' from (\ref{eqntimedecomp}).
%
%
\begin{lemma}
\label{lemb}
We have
\begin{eqnarray*}
&&\sum_{x \neq y \neq z} \bigl(\p_{x,z}[ \tau_1(y) \leq\Gamma_n,
\tau
_2(y) > \Gamma_n] -
\p_{x,z}[ \tau_2(y) \leq\Gamma_n, \tau_1(y) > \Gamma_n]
\bigr)\wt{\pi
}(z;x,y) \\
&&\qquad = \bigl( 1 + O(\Delta_n) \bigr) \sum_{x \neq y} \bigl( f_c(x,y)
- \ol{f}_c\bigr)^2 + O(e^{-\gamma c} (T_\mix)^2).
\end{eqnarray*}
\end{lemma}
\begin{pf}
Note that
\begin{eqnarray*}
&&\p_{x,z}[\tau_1(y) \leq\Gamma_n, \tau_2(y) > \Gamma_n] - \p
_{x,z}[\tau_2(y) \leq\Gamma_n, \tau_1(y) > \Gamma_n]\\
&&\qquad=
\p_x[\tau_1(y) \leq\Gamma_n] - \p_z[\tau_2(y) \leq\Gamma_n].
\end{eqnarray*}
Let $\delta_1(x,y,z) = O(e^{-\gamma c} \Upsilon_n + (g(x,z) + g(y,z) +
\Delta_n)(g(x,y) + \Delta_n))$ be the error term from Theorem \ref
{thmrnestimate} and $\delta_2(x,y,z) = O(e^{-\gamma c} \Upsilon_n +
\Delta_n(g(x,y) + g(x,z)))$ be the error term from Lemma \ref
{lemmixingcancel}. Then we can rewrite the summation in the
statement of the lemma as
\begin{eqnarray*}
&&\sum_{x \neq y \neq z}\bigl( \p_{x}[ \tau_1(y) \leq\Gamma_n] -
\p_{z}[ \tau_2(y) \leq\Gamma_n] \bigr)\wt{\pi}(z;x,y)\\
&&\qquad= \frac{1}{|V_n|}\sum_{x \neq y \neq z} \bigl( \p_x[\tau_1(y)
\leq
\Gamma_n] -
\p_z[ \tau_2(y) \leq\Gamma_n]\bigr) \bigl(1 +\varepsilon(x,y,z)\bigr)  \\
&&\qquad\quad{}+\frac{1}{|V_n|} \sum_{x \neq y \neq z} \bigl( |f_c(x,y) - f_c(y,z)| +
\delta_2(x,y,z)\bigr) \delta_1(x,y,z)\\
&&\qquad\equiv B_1 + B_2,
\end{eqnarray*}
where, by Theorem~\ref{thmrnestimate},
\[
\varepsilon(x,y,z) = \bigl(1 + O(\Delta_n)\bigr)\bigl( 2\ol{f}_c -
f_c(x,z) -
f_c(y,z)\bigr).
\]
Applying Assumption~\ref{assumpmain} repeatedly, it is tedious but
not difficult to see that $B_2 = O(e^{-\gamma c} T_\mix^2)$. By Lemma
\ref{lemmixingcancel},
\begin{eqnarray*}
B_1
&=& \bigl(1 + O(\Delta_n) \bigr)\frac{1}{|V_n|}\sum_{x \neq y \neq z}
\bigl( f_c(x,y) - f_c(y,z) + \delta_2(x,y,z) \bigr)\\
&&\hspace*{104.5pt}{}\times \bigl(2\ol{f}_c - f_c(x,z) - f_c(y,z)\bigr).
\end{eqnarray*}
Multiplying through, using the symmetry of $f$ in its arguments and
canceling many terms, this becomes
\[
\bigl( 1 + O(\Delta_n) \bigr) \sum_{x \neq y} \bigl( f_c(x,y) - \ol{f}_c\bigr)^2
+ O(e^{-\gamma c} (T_\mix)^2).
\]
\upqed\end{pf}

We will now show that the part of (\ref{eqnvarianceformula}) coming
from ``$A$'' of (\ref{eqntimedecomp}) is negligible. Roughly, the
reason for this is that it is unlikely for both walks to hit $y$
quickly, though in order to get a sufficiently good bound we will need
to take advantage of some more cancellation. This will in turn require
us to invoke (\ref{eqnradonsimplebound}), which is a rough estimate
of the Radon--Nikodym derivative of $\wt{\pi}$ with respect to $\pi$.
%
%
\begin{lemma}
\label{lema}
Uniformly in $n$, we have
\begin{eqnarray*}
&&\sum_{x \neq y \neq z} \bigl( \p_{x,z}[ \tau_1(y) < \tau_2(y)
\leq
\Gamma_n] -
\p_{x,z}[ \tau_2(y) < \tau_1(y) \leq\Gamma_n] \bigr) \wt{\pi
}(z;x,y)\\
&&\qquad = o(T_\mix^2).
\end{eqnarray*}
\end{lemma}
\begin{pf}
By (\ref{eqnradonsimplebound}), the summation in the statement of
the lemma is equal to
\begin{eqnarray*}
&&\frac{1}{|V_n|}\sum_{x \neq y \neq z} \bigl(\p_{x, z}[ \tau_1(y) <
\tau
_2(y) \leq\Gamma_n]
- \p_{x,z}[ \tau_2(y) < \tau_1(y) \leq\Gamma_n]\bigr)\\
&&\hspace*{46pt}{}\times\bigl(1+O\bigl(g(x,y) + g(y,z) + g(x,z)\bigr)\bigr).
\end{eqnarray*}
By symmetry, we see that this is equal to
\begin{eqnarray*}
&&\frac{1}{|V_n|}\sum_{x \neq y \neq z} \bigl(\p_{x, z}[ \tau_1(y) <
\tau
_2(y) \leq\Gamma_n]
- \p_{x,z}[ \tau_2(y) < \tau_1(y) \leq\Gamma_n]\bigr)\\
&&\hspace*{46pt}{}\times \bigl(O\bigl(g(x,y) + g(y,z) + g(x,z)\bigr)\bigr).
\end{eqnarray*}
Using
\[
\p_{x,z}[ \tau_1(y) < \tau_2(y) \leq\Gamma_n] \leq\p_{x,z}[\tau_1(y)
\leq\Gamma_n, \tau_2(y) \leq\Gamma_n]
\]
and the independence of $X_1, X_2$, we have the further bound
\[
\frac{1}{|V_n|}\sum_{x \neq y \neq z} O\bigl(\bigl[g(x,y) + \Delta_n\bigr)\bigl(g(y,z) +
\Delta_n\bigr)\bigr)O\bigl(g(x,y) + g(y,z) + g(x,z)\bigr).
\]
By the symmetry of $g$ in its arguments, we can rewrite this as
\begin{eqnarray*}
&&\frac{1}{|V_n|}\sum_{x \neq y \neq z} O\bigl( g(x,y) g^2(y,z) + g^2(x,y)
\Delta_n + g(x,y) \Delta_n^2\\
&&\qquad\hspace*{39.2pt}{}  + g(x,y)g(x,z) \Delta_n + g(x,y)g(x,z)g(y,z)\bigr).
\end{eqnarray*}
The terms in the summation are of order
\[
\CS_n T_\mix,\qquad {\CS_n T_\mix\log}|V_n|,\qquad
\frac{T_\mix^3 ({\log}|V_n|)^2}{|V_n|},\qquad
\frac{{T_\mix^3 \log}|V_n|}{|V_n|},\qquad
\frac{T_\mix^3}{|V_n|},
\]
respectively.
Assumption~\ref{assumpmain} implies that all of these are $o(T_\mix
^2)$, which gives the lemma.
\end{pf}

We complete this subsection by proving that ``$C$'' from (\ref
{eqntimedecomp}) is also negligible in comparison to the bound we
seek to prove. The intuition for this is that by time $\Gamma_n$, both
walks are very well mixed hence given that both have not hit $y$, the
difference in the probability that one hits before the other is of
smaller order than any negative power of $|V_n|$ (though we choose to
write $-100$). The proof will be in a slightly different spirit than
the previous lemmas.
%
%
\begin{lemma}
\label{lemc}
For any fixed $x,z$, we have
\[
\p_{x,z}[ \Gamma_n < \tau_1(y) < \tau_2(y)] -
\p_{x,z}[ \Gamma_n < \tau_2(y) < \tau_1(y)] = O(|V_n|^{-100}).
\]
\end{lemma}
\begin{pf}
We may assume without loss of generality that $x,z \neq y$.
The idea of the proof is to use a standard coupling argument to\vadjust{\goodbreak} show
that, conditional on $\{\tau_1(y) \wedge\tau_2(y) \geq\Gamma_n\}$,
the laws of $X_1(\Gamma_n)$ and $X_2(\Gamma_n)$ have total variation
distance $O(|V_n|^{-100})$ independent of $x,z$. To this end, we set
$\mu(z;x,y) = \p_x[ X(\Gamma_n) = z | \tau(y) \geq\Gamma_n]$.
Let $Y(t)$ be the process given by $X(t)$ conditioned on the event $\{
\tau(y) \geq\Gamma_n\}$. Then $Y(t)$ is Markov (though
time-inhomogeneous) as
\begin{eqnarray*}
&&\p[ Y(t) = z | Y(0) = z_0,\ldots,Y(t-1) = z_{t-1}]\\
&&\qquad= \p[ X(t) = z | X(0) = z_0,\ldots,X(t-1) = z_{t-1}, \tau(y) \geq
\Gamma_n]\\
&&\qquad= \frac{\p[ X(t) = z, \tau(y) \geq\Gamma_n| X(0) = z_0,\ldots,
X(t-1) = z_{t-1}]}{\p[\tau(y) \geq\Gamma_n | X(0) = z_0,\ldots, X(t-1)
= z_{t-1}]}\\
&&\qquad= \frac{\p_{z_{t-1}}[ X(1) = z, \tau(y) \geq\Gamma_n-(t-1)]}{\p
_{z_{t-1}}[\tau(y) \geq\Gamma_n-(t-1)]}
\end{eqnarray*}
depends only on $z,z_{t-1}$.
Recall that $T_k = k T_\mix$. For $t = c T_k$, note that
\begin{eqnarray*}
\nu(z;t,x)
&\equiv&\p_x[Y(t) = z]
= \frac{\p_x[X(t) = z, \tau(y) \geq\Gamma_n | \tau(y) \geq t]}{\p
_x[\tau(y) \geq\Gamma_n| \tau(y) \geq t]}\\
&=& \frac{\p_z[\tau(y) \geq\Gamma_n -t ] \p_x[ X(t) = z| \tau(y)
\geq
t]}{\p_x[\tau(y) \geq\Gamma_n| \tau(y) \geq t]}.
\end{eqnarray*}
Combining part~\ref{assumpmainhitmiss} of Assumption \ref
{assumpmain} with Lemma~\ref{lemconditionnothitrn}, we have that
\[
\frac{\p_z[\tau(y) \geq\Gamma_n -t ]}{\p_x[\tau(y) \geq\Gamma
_n| \tau
(y) \geq t]} = \Theta(1) \qquad\mbox{for } z \neq y.
\]

Also, since $\sum_{z} g(y,z) = T_\mix$ and $T_\mix= o(|V_n|)$, it
follows that for each $\varepsilon>0$ fixed, with $A = \{ z\dvtx g(y,z)
\leq\varepsilon\}$ we have $|A| / |V_n| = 1-o(1)$. Lemma~\ref
{lemconditionnothitrn} also implies that $\p_x[X(t) = z| \tau(y)
\geq t] = \Theta(1) \pi(z)$ on $A$ uniformly in $n$ large provided that
$k$ is large enough and $\varepsilon> 0$ is sufficiently small. This
implies that we can couple together the laws of $Y_u(c T_k), Y_v(c T_k)$
starting at $u,v$ distinct so that with probability $\rho> 0$, we have
$Y_u(c T_k) = Y_v(c T_k)$. If we iterate this procedure $c_1 = {\frac
{c_0}{\eta} \log}|V_n|$ times, $\eta= \eta(c,k,\rho)$, we get that
with probability $1 - O(|V_n|^{-c_1})$, we have $Y_u(\Gamma_n) =
Y_v(\Gamma_n)$. Consequently, we may assume that $c_0$ is sufficiently
large so that
\[
\max_{u,v} \| \nu(\cdot;\Gamma_n,u) - \nu(\cdot,;\Gamma_n,v) \|
_{\mathrm{TV}} =
O(|V_n|^{-500}).
\]
Let $\ol{\nu}$ be a measure so that $\max_{u} \| \nu(\cdot;\Gamma_n,u)
- \ol{\nu} \|_{\mathrm{TV}} = O(|V_n|^{-500})$.
Let $D = \{ \tau_1(y) \wedge\tau_2(y) \geq\Gamma_n\}$. Then we
have that
\begin{eqnarray*}
&&\p_{x,z}[ \Gamma_n < \tau_1(y) < \tau_2(y)] - \p_{x,z}[\Gamma_n
< \tau
_2(y) < \tau_1(y)]\\
&&\qquad= \bigl( \p_{x,z}[G_{12}(y) | D] - \p_{x,z}[G_{21}(y) | D] \bigr)
\p
_{x,z}[D]\\
&&\qquad= \sum_{u,v} \bigl(\p_{u,v}[ G_{12}(y)] - \p_{u,v} \p[ G_{21}(y)]
\bigr) \ol{\nu}(u;\Gamma_n,x) \ol{\nu}(v;\Gamma_n,z) \p_{x,z}[D]\\
&&\qquad\quad{} +
O(|V_n|^{-200})\\
&&\qquad= O(|V_n|^{-200}).
\end{eqnarray*}
\upqed\end{pf}
\begin{pf*}{Proof of Theorem~\ref{thmmain}}
To finish the proof of Theorem~\ref{thmmain}, we need to estimate
the diagonal (\ref{eqnvardiag}) and take care of the term in (\ref
{eqnvarianceneg}). Observe that (\ref{eqnvardiag}) is equal to
%
%
\begin{equation}\label{eqndiagest1}
2 \sum_{x} \p[G_{12}(x)] = |V_n| + O\biggl( \sum_x \p[\tau_1(x) =
\tau
_2(x)] \biggr).
\end{equation}
We can estimate the sum on the right-hand side using
\begin{eqnarray*}
\sum_x \p[\tau_1(x) = \tau_2(x)]
&=& \sum_x \sum_t \p[\tau_1(x) = t]\p[\tau_2(x) = t]\\
&\leq&\sum_x \sum_t \p[\tau_1(x) = t] \p[X_2(t) = x]
= 1.
\end{eqnarray*}
On the other hand, note
%
%
\begin{equation}\label{eqndiagest2}
\sum_x \bigl(f_c(x,x) - \ol{f}_c\bigr)^2 = \sum_x (1 - 2\ol{f}_c +
\ol{f}{}^2_c).
\end{equation}
By a union bound, we have $\ol{f}_c = O(c \Upsilon_n)$. Thus the
diagonal term in (\ref{eqndiagest1}) and (\ref{eqndiagest2})
differ by $O(c T_\mix) = o(e^{-\gamma c} T_\mix^2)$ (recall from
Assumption~\ref{assumpmain} that $T_\mix\to\infty$ as $n \to
\infty
$). This takes care of (\ref{eqnvardiag}).

We now turn to (\ref{eqnvarianceneg}). The previous lemma implies
\begin{eqnarray*}
&&\sum_{x \neq y \neq z} |\p_{x,z}[ G_{12}(y)] - \p
_{x,z}[G_{21}(y)]|\\
&&\qquad= \sum_{x \neq y \neq z} |\p_{x,z}[G_{12}(y), \tau_{1}(y) \leq
\Gamma
_n] - \p_{x,z}[ G_{21}(y), \tau_{2}(y) \leq\Gamma_n]|\\
&&\qquad\quad{} + O(|V_n|^{-50}).
\end{eqnarray*}
Observe $\{G_{ij}(y), \tau_i(y) \leq\Gamma_n\} \subseteq\{ \tau_i(y)
\wedge\tau_j(y) \leq\Gamma_n\} \cup\{ \tau_i(y) \leq\Gamma_n <
\tau
_j(y)\}$. Thus we can bound from above the previous expression by
\begin{eqnarray*}
&&\sum_{x \neq y \neq z}\bigl(2 \p_{x,z}[ \tau_1(y) \wedge\tau_2(y)
\leq\Gamma_n] +
|\p_{x}[ \tau_1(y) \leq\Gamma_n] - \p_{z}[ \tau_2(y) \leq\Gamma
_n]|\bigr)\\
&&\qquad\equiv E_1 + E_2.
\end{eqnarray*}
The term corresponding to $E_1$ can be bounded in a similar manner as
``$A$'' in the proof of Lemma~\ref{lema}. Indeed, by the independence
of $X_1,X_2$, we have that
\[
\p_{x,z}[ \tau_1(y) \wedge\tau_2(y) \leq\Gamma_n] \leq
\bigl(g(x,y) + \Delta_n\bigr) \bigl(g(y,z) + \Delta_n\bigr),
\]
which, when summed over $x,y,z$, is of order $O( ({\log}|V_n|)^2 |V_n|
T_\mix^2)$. We can estimate $E_2$ using techniques similar to the proof
of Lemma~\ref{lemb} since by Lemma~\ref{lemmixingcancel},
\[
| \p_x[ \tau_1(y) \leq\Gamma_n] - \p_z[ \tau_2(y) \leq\Gamma
_n]| =
O\bigl(g(x,y) + g(y,z) + \delta_2(x,y,z)\bigr),
\]
where, as in the proof of Lemma~\ref{lemb}, $\delta_2(x,y,z)$
corresponds to the error from Lemma~\ref{lemmixingcancel}. When
summed over $x,y,z$, this is of order $O(|V_n|^2 T_\mix)$.
Therefore
\[
(E_1 + E_2) \biggl( o \biggl(\frac{T_\mix}{{|V_n|^2 \log}|V_n|}
\biggr) +
O\biggl(\frac{1}{|V_n|^2} \biggr)\biggr) = o(T_\mix^2)
\]
as desired.
\end{pf*}

\section{Further questions}
\label{secproblems}

\begin{figure}[b]
\begin{tabular}{@{}c@{\hspace*{4pt}}c@{\hspace*{4pt}}c@{}}

\includegraphics{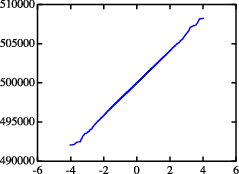}
 & \includegraphics{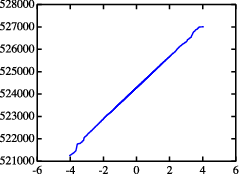} & \includegraphics{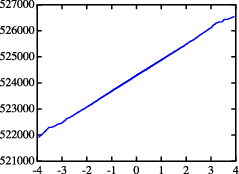}\\
(a) & (b) & (c)
\end{tabular}
\caption{$Q$--$Q$ plots based on 20,000 simulations of the number of
sites visited by $X_1$ before $X_2$ against an appropriately fitted
normal distribution, supporting the conjecture of asymptotic
normality, where \textup{(a)} $\Z_{100}^3$; \textup{(b)} $\Z_{32}^4$
and \textup{(c)} $\Z_2^{20}$.} \label{figqqplots}
\end{figure}

(1) The first step in proving a sequence of random variables
$(X_n)$ has a Gaussian limit after appropriate normalization is the
determination of the asymptotic mean and variance. We remarked in the
beginning that, in our case, the expected number of sites painted 1
is $|V_n|/2$, and Theorem~\ref{thmmain} gives the limiting variance.
Figure~\ref{figqqplots} shows $Q$--$Q$ plots of the empirical
%
%
distribution of the number of sites painted $1$ in the final coloring
against an appropriately fitted normal for three different base graphs.
Based on these plots, we conjecture that
\[
\frac{ |\CA_i| - \E|\CA_i|}{\sqrt{\var(|\CA_i|)}}
\]
has a normal limit for all graphs satisfying Assumption
\ref{assumpmain}.

(2) Our derivation of the variance ignores the time aspect of
the problem in the sense that it gives no indication of at what point
in the process of coverage the variance is ``created.'' Does it come in
bursts or continuously? Does it come sooner than any multiple of the
cover time or perhaps in $[\varepsilon T_\cov,T_\cov]$? More generally,
when normalized appropriately, does the the process $t \mapsto\sum_x
\one_{\{\tau_1(x) < \tau_2(x) \leq t\}}$ have a scaling limit?\vadjust{\goodbreak}

(3) We make repeated used of the symmetry afforded by the fact
that we consider two random walks moving at the same speed on vertex
transitive graph. It would be interesting to see if a similar result
holds when the various degrees of symmetry are broken. Starting points
for exploring this problem include considering continuous time walks
moving at various speeds, multiple walks and graphs which are not
vertex transitive.

(4) Theorem~\ref{thmtori} only holds for tori of dimension $d
\geq3$ as the case $d=2$ falls just outside of the scope of Theorem
\ref{thmmain}. It would be interesting to see a more refined analysis
carried out to handle this case.

(5) That the variance computed in Theorem~\ref{thmtori} for
$d=3,4$ is significantly larger than in the i.i.d. case suggests that
the clusters which have an unusually large number of sites painted a
given color are either larger or more dense than in an i.i.d. marking.
How large and frequent are such clusters? What is their geometric structure?

(6) Another interesting quantity is the size $\CB$ of the
boundary separating the sites painted $1$ and $2$, as studied in
\cite{JLSH96}. It is not difficult to see that there exists a constant
$\beta
_d > 0$ such that $\E|\CB| \sim\beta_d n^d$ when $d \geq3$ as $n
\to
\infty$. Indeed, this follows since the probability that $\{\tau_1(y) <
\tau_2(y)\}$ for $y \sim x$ given $\{\tau_1(x) < \tau_2(x)\}$ converges
to a limit $p_d \in(0,1)$. Note that this is of the same order of
magnitude as $\E| \CA_1|$. Is it also true that $\var(|\CB|) =
\Theta
(\var(|\CA_1|))$ or do these quantities differ significantly?

\section*{Acknowledgments}

I am grateful to Yuval Peres for introducing me to this problem as well
as for helpful discussions on different methods of proof over the
course of a summer internship in the Theory Group at Microsoft
Research. The idea of considering the symmetrized problem described in
the beginning of Section~\ref{secvar} was crucial and motivated the use
of symmetrization whenever possible elsewhere in the proof. In
addition, I thank Robin Pemantle for communicating the conjecture that
led to Theorem~\ref{thmtori}, which was formulated jointly with
Peres~\cite{PP09}. Having such an accurate target was rather helpful
for crafting the proof, in particular for achieving a precise lower
bound. I thank Amir Dembo for pointing out that the original proof of
Theorem~\ref{thmtori} gave both the constant in front of the leading
term and held in much greater generality, which yielded Theorem
\ref{thmmain}. Finally, I~thank both Amir Dembo and Perla Sousi for
many useful comments on an earlier draft of this article.


%

\printaddresses

\end{document}